\documentclass[journal]{IEEEtran}
\ifCLASSINFOpdf
\else
 \fi

\usepackage[usenames]{color}
\usepackage[colorlinks = true]{hyperref}
\usepackage[english]{babel}
\usepackage{bm,bbm,amsmath,amssymb,amsthm,graphicx}
\usepackage{enumitem}
\usepackage{url,cite}
\usepackage{pifont}
\usepackage[caption=false,font=normalsize,labelfont=sf,textfont=sf]{subfig}
\usepackage{algorithm}
\usepackage{algpseudocode}
\usepackage{tikz}
\usetikzlibrary{calc,patterns,decorations.pathmorphing,decorations.markings}

\newcommand{\SLb}[1]{\textcolor{black}{#1}}
\newcommand{\SLbb}[1]{\textcolor{black}{#1}}
\newcommand{\MBWbb}[1]{\textcolor{black}{#1}}
\newcommand  \stack[2]   {\overset{\text{#1}}{#2}}

\newcommand{\change}[1]{\ensuremath{\operatorname{#1}}}
\newcommand{\MAT}{\left[ \begin{array}}
\newcommand{\mat}{\end{array} \right]}
\newtheorem{Definition}{Definition}[section]
\newtheorem{Corollary}{Corollary}[section]
\newtheorem{Lemma}{Lemma}[section]
\newtheorem{Theorem}{Theorem}[section]

\newtheorem{Q}{Question}[section]
\newtheorem{Proposition}{Proposition}
\newtheorem{Remark}{Remark}[section]

\def \st {\operatorname*{s.t. }}

\def \a {\bm{a}}
\def \A {\mathbf{A}}
\def \AA {\mathcal{A}}

\def \b {\bm{b}}
\def \B {\mathbf{B}}

\def \CC {\mathcal{C}}

\def \D {\mathbf{D}}
\def \DD{\mathcal{D}}

\def \e {\bm{e}}

\def \E {\mathbf{E}}

\def \EEE{\mathbb{E}}

\def \hwt {\widetilde{h}}

\def \I {\mathbf{I}}

\def \K {\mathbf{K}}

\def \M {\mathbf{M}}

\def \N {\mathbf{N}}
\def \NN {\mathcal{N}}

\def \PP {\mathcal{P}}

\def \RRR {\mathbb{R}}
\def \S {\mathbf{S}}

\def \SSS {\mathbb{S}}
\def \St {\widetilde{\mathbf{S}}}

\def \TT {\mathcal{T}}

\def \u {\bm{u}}

\def \ut {\widetilde{\bm{u}}}
\def \U {\mathbf{U}}

\def \Ut {\widetilde{\mathbf{U}}}
\def \v {\bm{v}}
\def \V {\mathbf{V}}

\def \x {\bm{x}}
\def \xs {\bm{x}^\star}

\def \X {\mathbf{X}}

\def \Xh {\widehat{\mathbf{X}}}
\def \Xs {\mathbf{X}^\star}

\def \y {\bm{y}}
\def \Y {\mathbf{Y}}

\def \sumn{\sum_{n=1}^N}

\def \tr{\change{tr}}
\def \diag{\change{diag}}
\def \St{\change{St}}

\def \bPhi {\boldsymbol{\Phi}}

\def \bxi {\boldsymbol{\xi}}

\def \bLambda {\boldsymbol{\Lambda}}

\def \zero {\mathbf{0}}

\begin{document}

\title{Landscape Correspondence of Empirical and Population Risks in the Eigendecomposition Problem}

\author{Shuang~Li,
        Gongguo Tang,
        and Michael B. Wakin
\thanks{SL is with the Department of Mathematics, University of California, Los Angeles, CA 90095. Email: shuangli@math.ucla.edu. GT and MBW are with the Department of Electrical Engineering, Colorado School of Mines, Golden, CO 80401. Email: \{gtang,mwakin\}@mines.edu.}
}

\maketitle

\begin{abstract}
Spectral methods include a family of algorithms related to the eigenvectors of certain data-generated matrices.
  In this work, we are interested in studying the geometric landscape of the eigendecomposition problem in various spectral methods. In particular, we first extend known results regarding the landscape at critical points to larger regions near the critical points in a special case of finding the leading eigenvector of a symmetric matrix. For a more general eigendecomposition problem, inspired by recent findings on the connection between the landscapes of empirical risk and population risk, we then build a novel connection between the landscape of an eigendecomposition problem that uses random measurements and the one that uses the true data matrix. We also apply our theory to a variety of low-rank matrix optimization problems and conduct a series of simulations to illustrate our theoretical findings.


\end{abstract}


\IEEEpeerreviewmaketitle


\section{Introduction}
\label{sec:intr}


\IEEEPARstart{S}{pectral} methods are of fundamental importance in signal processing and machine learning due to their simplicity and effectiveness. They have been widely used to extract useful information from noisy and partially observed data in a variety of applications, including dimensionality reduction~\cite{saul2006spectral}, tensor estimation~\cite{cai2019nonconvex},
ranking from pairwise comparisons~\cite{negahban2017rank}, low-rank matrix estimation~\cite{keshavan2010matrix}, and community detection~\cite{newman2006finding}, to name a few.
As is well known, the theoretical performance of the gradient descent algorithm and its variants is heavily dependent on a proper initialization~\cite{sutskever2013importance,ruder2016overview}. In recent years, spectral methods and their variants are commonly used as an initialization step for many algorithms in order to guarantee linear convergence~\cite{candes2015phase,ma2019implicit,wang2017solving,chi2019nonconvex}. Spectral initialization approaches have proven to be very powerful in providing a ``warm start'' for many non-convex matrix optimization problems such as matrix sensing~\cite{chi2019nonconvex,li2019non}, matrix completion~\cite{sun2016guaranteed,chen2019nonconvex}, blind deconvolution~\cite{li2019rapid,lee2018fast}, phase retrieval~\cite{candes2015phase,bendory2017non}, quadratic sensing~\cite{sanghavi2017local,li2019nonconvex}, joint alignment from pairwise differences~\cite{chen2018projected}, and so on. In these problems, spectral methods provide an initialization by computing the leading eigenvector(s) of a certain surrogate matrix constructed from the given measurements.


In signal processing and machine learning problems involving randomized data, the empirical risk is a cost function that is typically written as a sum of losses dependent on each sample of the true data~\cite{vapnik2006estimation,koltchinskii2001rademacher,vapnik1992principles}. The population risk is then defined as the expectation of the empirical risk with respect to the randomly generated data. Since there is no randomness in the population risk, it is often easier to analyze. The correspondence between the geometric landscape of empirical risk and that of its corresponding population risk was first studied in~\cite{mei2018landscape} for the case when the population risk is strongly Morse, that is, when the Hessian of the population risk has no zero eigenvalues at or near the critical points. A subsequent work~\cite{li2019landscape} extends this to a more general framework by removing the strongly Morse assumption. In particular, the uniform convergence of the empirical risk to the population risk is studied in a situation where the population risk is degenerate, namely, the Hessian of the population risk can have zero eigenvalues.

Inspired by these recent findings, in this work, we are interested in analyzing the landscape of the eigendecomposition problem widely used in the spectral methods. \SLb{In particular, we formulate the eigendecomposition problem as an optimization problem constrained on the Stiefel manifold.}
We view the quadratic eigendecomposition cost function based on the true symmetric data matrix as a population risk since it is deterministic. Accordingly, the cost function based on random measurements \SLb{(via the surrogate matrix)} is then viewed as the empirical risk. Under certain assumptions, we can guarantee that the population risk is the expectation of the empirical risk.
\SLb{The problem of finding the  leading eigenvector has been studied in the literature, and the landscape of the population risk {\em at} the critical points can be found in~\cite[\SLbb{Section 4.6.2}]{absil2009optimization}. Here, we extend this geometric analysis {\em at} the critical points to larger regions {\em near} the critical points. We also extend these results to a more general eigendecomposition problem, in which we consider the problem of finding the $r$ leading eigenvectors. We again establish the favorable geometry of the population risk in this more general problem in regions both {\em at} and {\em near} the critical points.}
We then build a connection between the critical points of the empirical risk and population risk, which makes it possible to directly study the landscape of the empirical risk using the landscape of the population risk.
To further support our theory, we apply it to a variety of low-rank matrix optimization problems such as matrix sensing, matrix completion, phase retrieval, and quadratic sensing.

A number of existing works conduct landscape analysis on a specific problem of interest~\cite{de2015global,ge2015escaping,ge2017no,sun2018geometric,li2020global,li2019geometry}. However, we focus more on the general eigendecomposition problems that are commonly used in the spectral methods that can be applied to solve these problems. We believe that we are the first to build this connection between the landscape of the empirical and population eigendecomposition problems.
\SLb{This can be of its own independent interest.}


The remainder of this paper is organized as follows. In Section~\ref{sec:prob}, we formulate our problem by introducing the eigendecomposition problem.
In Section~\ref{sec:main}, we present our main results on the landscape of the population risk in the eigendecomposition problem and illustrate how to infer the landscape of an empirical risk from its  population risk.
We then survey a wide range of applications in Section~\ref{sec:appl} and conduct a series of experiments in Section~\ref{sec:simu} to further support our theoretical analysis.
Finally, we conclude our work and discuss future directions in Section~\ref{sec:conc}.

\textbf{Notation:}
We use boldface uppercase letters (e.g., $\X$) and boldface lowercase letters (e.g., $\x$) to denote matrices and vectors, respectively. Scalars or entries of matrices and vectors are not bold.
For example, we let $\{\x_s = \X(:,s)\}_{s=1}^r$ denote the columns of a matrix $\X\in\RRR^{N\times r}$ and $\x_s(j)$ denote the $j$-th entry of $\x_s \in \RRR^{N}$, i.e., $\x_s(j) = \X(j,s)$ denotes the $(j,s)$-th entry of $\X$. Similarly, we also use $\X(j,:)$ to denote the $j$-th row of $\X$.
For any matrix $\X$, let $\|\X\|_F$, $\|\X\|$, and $\|\X\|_\infty$ denote the Frobenius norm, spectral norm, and maximum entry (in absolute value) of $\X$, respectively. \SLbb{We use $\langle \A,\B \rangle = \tr(\B^\top \A)$ and $\langle \a,\b \rangle = \b^\top \a$ to denote the  standard inner product for matrices and vectors in Euclidean space, respectively.} We use $C$, $c$, $c_1$, $c_2,\ldots$ to denote numerical constants with values that may change from line to line.
The matrix commutator is defined as $[\X_1,\X_2] \triangleq \X_1 \X_2-\X_2 \X_1$. Denote $[N] \triangleq \{1,\cdots,N\}$.
\MBWbb{We use $\PP_{mu}(\cdot)$ to indicate permutations; for example, $\{1,2,3\} = \PP_{mu}(\{1,2,3\})$ and $\{2,1,3\} = \PP_{mu}(\{1,2,3\})$, but $\{2,1,4\}  \neq \PP_{mu}(\{1,2,3\})$.}
Finally, we denote $\e_i \in \RRR^{N}$ as the $i$-th column of an identity matrix $\I_N$ and $\mathcal{I}$ as the identity operator.

\section{Problem Formulation}
\label{sec:prob}

In this work, we consider the following eigendecomposition problem
\begin{align}
\max_{\X\in\RRR^{N\times r}}~\tr(\X^\top \M \X \N	) \quad \st~\X^\top\X = \I_r,
\label{eqn:max_eigen}
\end{align}
where $\M\in \RRR^{N\times N}$ is a symmetric matrix with eigenvalues satisfying $\lambda_1 > \lambda_2 >\cdots > \lambda_r > \lambda_{r+1} \geq \cdots \geq \lambda_N$.\footnote{\SLb{When $\M$ is non-symmetric, one can show that maximizing $\tr(\X^\top \M \X \N)$ is equivalent to maximizing $\tr( \X^\top(\M+\M^\top)\X\N)$. Since $(\M+\M^\top)$ is symmetric, the maximizer of $\tr(\X^\top \M \X \N)$ therefore corresponds to the leading eigenvectors of $(\M+\M^\top)$ instead of the original $\M$.}}
$\N = \diag([\mu_1~\mu_2~\cdots~\mu_r])\in \RRR^{r \times r}$ is a diagonal matrix. A recent paper~\cite{birtea2019first} shows that if the diagonal entries in $\N$ are pairwise distinct and strictly positive, then $\X$ is a critical point if and only if the columns of $\X$ are eigenvectors of $\M$. This inspires us to set $\N$ as $\N = \diag([r~r-1~\cdots~1])$, i.e., $\mu_i = r-i+1$ for all $i\in[r]$.
Observe that the maximization problem~\eqref{eqn:max_eigen} is equivalent to the following minimization problem
\begin{align}
\min_{\X\in\RRR^{N\times r}}~-\frac 1 2 \tr(\X^\top \M \X \N	)	\quad \st~\X^\top\X = \I_r,
\label{eqn:min_eigen}
\end{align}
which can be viewed as minimizing a population risk
\begin{align}
g(\X) = -\frac 1 2 \tr(\X^\top \M \X \N	)
\label{eqn:def_pr}
\end{align}
on the Stiefel manifold $\St(N,r) \triangleq \{\X\in\RRR^{N\times r}: \X^\top \X=\I_r\}$. The above population risk~\eqref{eqn:def_pr} is also known as the Brockett cost function. In future sections this deterministic quantity will arise as the expectation of an empirical risk in certain learning problems.

\section{Main Results}
\label{sec:main}

\subsection{Warm Up: Special Case with $r=1$}

We begin with the simple example of finding the leading eigenvector,
namely
\begin{align}
\min_{\x\in\RRR^{N}}~-\frac 1 2 \x^\top \M \x 	\quad \st~\|\x\|_2 = 1,
\label{eqn:min_eigen_r1}
\end{align}
which is a special case of the eigendecomposition problem~\eqref{eqn:min_eigen} with $r=1$. This problem can be viewed as minimizing a population risk
\begin{align}
g(\x) = -\frac 1 2\x^\top \M \x
\label{eqn:def_pr_r1}
\end{align}
on the unit sphere $\SSS^{N-1} \triangleq \{\x\in\RRR^N: \|\x\|_2=1\}$. Since analysis on the unit sphere is often easier than analysis on the Stiefel manifold, we focus on this special case first.

Problem~\eqref{eqn:min_eigen_r1} is widely used in spectral methods to compute the leading eigenvector or to generate a reasonably good initialization for  low-rank matrix optimization problems~\cite{chi2019nonconvex}.
We summarize the characterization of the critical points of the population risk~\eqref{eqn:def_pr_r1} in the following proposition.
\begin{Proposition}\label{prop:eigvec}
\cite[Propositions 4.6.1, 4.6.2]{absil2009optimization}
For a symmetric matrix $\M\in\RRR^{N\times N}$, a vector $\x\in\RRR^{N}$ with $\|\x\|_2 = 1$ is an eigenvector of $\M$ if and only if it is a critical point of the population risk~\eqref{eqn:def_pr_r1}\footnote{Here, the population risk~\eqref{eqn:def_pr_r1} is constrained on the unit sphere.}. Moreover, denote $\{\lambda_n\}_{n=1}^N$ with $\lambda_1 > \lambda_2 \geq \lambda_3 \geq \cdots \geq \lambda_N$ as the eigenvalues of $\M$ and $\{\v_n\}_{n=1}^N$ as the associated eigenvectors. Then, we have
\begin{itemize}[leftmargin=*]
\item $\pm \v_1$ are the only global minimizers of~\eqref{eqn:def_pr_r1}.
\item $\pm \v_N$ are global maximizers of~\eqref{eqn:def_pr_r1}. If $\lambda_N < \lambda_{N-1}$, then $\pm \v_N$ are the only global maximizers.
\item $\pm \v_n$ with $\lambda_N< \lambda_n<\lambda_1$ are  saddle points of~\eqref{eqn:def_pr_r1}.
\end{itemize}

\end{Proposition}

Note that the above proposition only characterizes the landscape of the population risk~\eqref{eqn:def_pr_r1} at the critical points.
Next, we extend this analysis to larger regions near the critical points, which further allows us to build a connection between the landscape of the population risk and the corresponding empirical risk in various applications involving random measurements of~$\M$.
With some elementary calculations, we can write the {\em Euclidean gradient and Hessian} of $g(\x)$ as
$
\nabla g(\x) = -\M \x,~ \text{and}~
\nabla^2 g(\x) = -\M.	$
The {\em Riemannian gradient and Hessian} on the unit sphere are then obtained from the projection of the Euclidean gradient and Hessian, i.e.,\footnote{\MBWbb{One can refer to \cite[Sections 4.6.1 and 5]{absil2009optimization} for more details on how to compute Riemannian gradients and Hessians.}}
\begin{align*}
\text{grad}~g(\x) & \!=\! (\I_N-\x\x^\top) \nabla g(\x) = (\x\x^\top - \I_N)\M\x,	\\
\text{hess}~g(\x)& \!=\! (\I_N\!\!-\!\x\x\!^\top\!)(\nabla^2 \!g(\x) \!-\! \x\!^\top \nabla \!g(\x) \I_N\!)(\I_N\!-\!\x\x\!^\top\!)\\
 &\!=\! (\I_N\!\!-\!\x\x\!^\top\!)(\x\!^\top\M\x\I_N\!\!-\!\M )(\I_N\!\!-\!\x\x\!^\top\!),
\end{align*}
 where $\I_N$ denotes the $N\times N$ identity matrix.

The following theorem, which establishes the favorable geometry of the population risk~\eqref{eqn:def_pr_r1} in regions near the critical points, is proved in Appendix~\ref{sec:proof_main_r1}.
\begin{Theorem}\label{thm:main_r1}
Assume that $\M\in \RRR^{N\times N}$ is a symmetric matrix with eigenvalues satisfying $\lambda_1 > \lambda_2 \geq \lambda_3 \geq \cdots \geq \lambda_N$. For the population risk $g(\x)$ defined in~\eqref{eqn:def_pr_r1}, there exist two positive numbers $\epsilon = 0.2(\lambda_1-\lambda_2)$ and $\eta = 0.3(\lambda_1-\lambda_2)$ such that
\begin{align*}
|\lambda_{\min}(\text{{\em hess}}~g(\x))| \geq \eta	
\end{align*}
holds if $\|\text{{\em grad}}~g(\x)\|_2 \leq \epsilon$.	
\end{Theorem}




Note that Proposition~\ref{prop:eigvec} indicates that the population risk~\eqref{eqn:def_pr_r1} satisfies the {\em strict saddle property}\footnote{As defined in~\cite{ge2015escaping}, a twice differentiable function satisfies the strict saddle property if the Hessian is positive definite when evaluated at any local minimum and contains a negative eigenvalue at any other stationary point.}, which allows many iterative algorithms to avoid saddle points and converge to local minima~\cite{ge2015escaping,du2017gradient,anandkumar2016efficient,jin2018accelerated,reddi2018generic}. Theorem~\ref{thm:main_r1} extends the results of Proposition~\ref{prop:eigvec} that hold at exact critical points to larger regions near critical points. This further establishes that the population risk~\eqref{eqn:def_pr_r1} satisfies the {\em robust strict saddle property} and guarantees that many local search algorithms can in fact converge to local minima in polynomial time~\cite{ge2015escaping,sun2015nonconvex,zhu2021global}. In summary, the population risk~\eqref{eqn:def_pr_r1} has a favorable geometry.


\subsection{General Case with $r\geq 1$}
With some elementary calculations \MBWbb{involving a Taylor expansion}, we can write the {\em Euclidean gradient and Hessian} of the population risk $g(\X)$ in~\eqref{eqn:def_pr} as
$\nabla g(\X) = -\M \X\N, ~ \nabla^2 g(\X) [\U,\U] = -\langle\M,\U\N\U^\top \rangle.$
The {\em Riemannian gradient and Hessian} on the Stiefel manifold are then obtained from the projection of the Euclidean gradient and Hessian~\cite{absil2009optimization}:
\begin{align*}
&\text{grad}~g(\X) = (\X\X^\top - \I_N)\M\X \N	-  \frac 1 2 \X[\X^\top \M \X, \N], \\
&\text{hess}~g(\X)[\U,\U] =  \langle \X^\top \M \X, \U^\top \U \N\rangle \!-\!\langle \M,\!\U\N\U\!^\top \!\rangle,
\end{align*}
 where $\U\in\RRR^{N\times r}$ is a matrix that belongs to the tangent space of $\St(N,r)$ at $\X$, i.e.,
\begin{align*}
\U \!\in \!\TT_{\X}&\St(N,r) =\! \{\U \in \RRR^{N\times r}: \X^\top \U + \U^\top \X \!= \zero\} \\
&~=\! \{\X \S+ \X_\perp \K: \S^\top \!=\! -\S, \K \in \RRR^{(N-r) \times r} \}.
\end{align*}
\SLbb{Here, $\X_{\perp}$ is any $N\times (N-r)$ matrix such that span$(\X_{\perp})$ is the orthogonal complement of span$(\X)$.}

It follows from~\cite[Section 4.8.2]{absil2009optimization} that a matrix $\X\in\RRR^{N\times r}$ with $\X^\top\X = \I_r$ is a critical point of the constrained population risk~\eqref{eqn:def_pr} if and only if its columns are eigenvectors of $\M$. The following theorems, which establish the favorable geometry of the population risk~\eqref{eqn:def_pr} with $r\geq 1$ in regions at and near the critical points, are proved in Appendices~\ref{sec:proof_main0} and~\ref{sec:proof_main}, respectively. Again, the  favorable geometry means that many iterative algorithms can avoid saddle points and converge to local minima in polynomial time, as in the case with $r=1$.

\begin{Theorem}\label{thm:main0}
Assume that $\M \in \RRR^{N\times N}$  with eigenvalues satisfying $\lambda_1 > \lambda_2 >\cdots > \lambda_r > \lambda_{r+1} \geq \cdots \geq \lambda_N $. Define $\Omega \triangleq \{i_1,\cdots i_r\}$ as a subset of $[N]$. Denote $\X_\Omega = [\x_{i_1},\cdots,\x_{i_r}] \in \RRR^{N\times r}$ with $\{\x_{i_j}\}_{j=1}^r$ being the $i_j$-th eigenvector of $\M$. Then, $\X_\Omega$ is a critical point of $g(\X)$ defined in~\eqref{eqn:def_pr}\footnote{Here, the population risk~\eqref{eqn:def_pr} is constrained on the Stiefel manifold.}.
Moreover, we have
\begin{itemize}[leftmargin=*]
\item $\Omega = [r]$: $\X_\Omega $ is a global minimizer of~\eqref{eqn:def_pr}. 
\item $\Omega = \PP_{mu}([r]) \neq [r]$: $\X_\Omega $ is a strict saddle point of~\eqref{eqn:def_pr}.
\item $\Omega \neq \PP_{mu}([r])$: $\X_\Omega $ is a strict saddle point of~\eqref{eqn:def_pr}.
\end{itemize}	
\end{Theorem}

\begin{Theorem}\label{thm:main}
Assume that $\M  \in \RRR^{N\times N}$  with eigenvalues satisfying $\lambda_1 > \lambda_2 >\cdots > \lambda_r > \lambda_{r+1} \geq \cdots \geq \lambda_N $. Define $d_{\min} \triangleq \min_{1\leq s<j \leq r+1} (\lambda_s-\lambda_j)$ as the minimal distance between any two of the first $r+1$ eigenvalues. For the population risk $g(\X)$ defined in~\eqref{eqn:def_pr}, there exist two positive numbers  $\epsilon = \frac{1}{72} r^{-1} d_{\min}$ and $\eta = 0.11 d_{\min}$ such that
\begin{align*}
|\lambda_{\min}(\text{{\em hess}}~g(\X))| \geq \eta	
\end{align*}
holds if $\|\text{{\em grad}}~g(\X)\|_F \leq \epsilon$.	
\end{Theorem}

\subsection{Landscape Correspondence between Empirical and Population Risks}
In many low-rank matrix optimization problems,
we may only have random measurements of the true matrix $\M$; several such applications are surveyed and detailed in Section~\ref{sec:appl}. In such problems, one may be able to construct a surrogate matrix $\Y$ from the measurements. This results in a corresponding empirical risk
\begin{align}
f(\X) = -\frac 1 2 \tr(\X^\top \Y \X \N	).
\label{eqn:def_er}
\end{align}
By taking expectation of $f(\X)$ with respect to the random measure used to obtain $\Y$, one often has $g(\X) = \EEE f(\X)$ \SLb{due to the fact that $\M = \EEE \Y$}. Namely, the population risk~\eqref{eqn:def_pr} is the expectation of the above empirical risk.

According to~\cite[Theorem 2.1, Corollary 2.1]{li2019landscape}, we can utilize Theorem~\ref{thm:main} above to build a connection between the critical points of the empirical risk~\eqref{eqn:def_er} and the population risk~\eqref{eqn:def_pr}.
We summarize this result in the following corollary.
\SLbb{Note that this corollary can be viewed as a direct result obtained by applying the more general theory in \cite[Theorem 2.1, Corollary 2.1]{li2019landscape} to our specific eigendecomposition problem.}
\begin{Corollary}\label{coro:main}
Denote $\DD$ as a maximal connected  and compact subset of the set $\{\X\in\St(N,r): \|\text{{\em grad}}~g(\X)\|_F \leq\epsilon\}$ with a $\CC^2$ boundary $\partial \DD$.\footnote{According to the proof of \cite[Theorem 2.1]{li2019landscape}, $\{\mathbf{X}\in \St(N,r):\|\text{grad}~ g(\mathbf{X})\|_F \leq \epsilon\}$ can be partitioned into disjoint connected compact components with each containing at most one local minimum. Here, $\mathcal{D}$ is one such component, so we can guarantee the connectedness and smoothness of the boundary for $\mathcal{D}$.}
Under the same assumptions used in Theorem~\ref{thm:main} and $g(\X) = \EEE f(\X)$, if the Riemannian gradient and Hessian of the population risk~\eqref{eqn:def_pr} and empirical risk~\eqref{eqn:def_er} satisfy\footnote{As introduced in Section~\ref{sec:appl}, these assumptions  are satisfied with high probability for suitable choices of $\epsilon$ and $\eta$. }
\begin{align}
&\sup_{\X\in \St(N,r)} \|\text{{\em grad}}~f(\X) - \text{{\em grad}}~g(\X)\|_F \leq \frac{\epsilon}{2},\label{eqn:as_grad}\\
&\sup_{\X\in \St(N,r)} \|\text{{\em hess}}~f(\X) - \text{{\em hess}}~g(\X)\| \leq \frac{\eta}{2},	\label{eqn:as_hess}
\end{align}
the following statements hold:
\begin{itemize}[leftmargin=*]
\item If $g$ has no local minima in $\DD$, then $f$ has no local minima in $\DD$.
\item If $g$ has one local minimum in $\DD$, then $f$ has one local minimum in $\DD$.
\MBWbb{Let $\X$ and $\Xh$ denote the local minima of the population risk~\eqref{eqn:def_pr} and the empirical risk~\eqref{eqn:def_er}, respectively, and suppose $\X \in \DD$ and $\Xh \in \DD$. Suppose the pre-image of $\DD$ under the exponential mapping $\text{Exp}_{\X}(\cdot)$ is contained in a ball at the origin of the tangent space $\TT_{\X}\St(N,r)$ with radius $\rho$.
Consider the differential of the exponential mapping $\text{DExp}_{\X}(\V)$, and let $\sigma$ be an upper bound on the operator norm of this differential for all $\V \in \TT_{\X}\St(N,r)$ with Frobenius norm less than $\rho$. Let $L_H$ be an upper bound on the Lipschitz constant of the pullback Hessian for the population risk at the origin of $\TT_{\X}\St(N,r)$.
Then, as long as $\epsilon \leq \frac{\eta^2}{2\sigma L_H}$, the Riemannian distance between the two local minima is upper bounded by $2\sigma(\rho)\epsilon/\eta$.}
\item If $g$ has strict saddle points in $\DD$, then if $f$ has any critical points in $\DD$, they must be strict saddle points.
\end{itemize}
\end{Corollary}

\begin{Remark}
	\textcolor{black}{Corollary~\ref{coro:main} depends on several parameters which we discuss here. The radius $\rho$ depends on $\epsilon$ (via $\DD$) and can therefore be made smaller by choosing smaller $\epsilon$; it has an upper bound of $0.89\pi$.\footnote{\MBWbb{We use the fact that the injectivity radius of the Stiefel manifold is at least $0.89\pi$~\cite{sutti2020leapfrog,rentmeesters2013algorithms}.}} Next, it can be shown that $\sigma \leq \sigma(\rho) := \exp(\rho^2 + 3\rho + 1) (2\rho^2 + 5\rho + 3)$.
	Finally, the parameter $L_H$ can be explicitly computed using ideas in the proof of \cite[Proposition 4.1]{zhang2018cubic}, which is independent of $\X$ and the dimension of the Stiefel manifold. }
\end{Remark}



To illustrate the result in Corollary~\ref{coro:main}, we introduce in the next section a series of applications where assumptions~\eqref{eqn:as_grad} and~\eqref{eqn:as_hess} can be satisfied.

\MBWbb{Before proceeding, we comment on the contributions of our main results. First, Theorems~\ref{thm:main0} and~\ref{thm:main} have relevance to the fundamental problem of computing the eigenvectors of a matrix $\M$. Theorem~\ref{thm:main0} extends the results from~\cite[Propositions 4.6.1, 4.6.2]{absil2009optimization} to the general case $r \ge 1$, and Theorem~\ref{thm:main} further establishes the robust strict saddle property for the landscape of the objective function $g(\X)$. As we have noted, this favorable geometry means that many iterative algorithms can avoid saddle points and converge to local minima in polynomial time. Meanwhile, classical methods such as the power method (when $r=1$) and subspace iteration method (when $r > 1$) also exist for computing eigenvectors of a matrix. Though it is beyond the scope of this paper, there are similarities to the power method and the subspace iteration method when deriving a Riemannian gradient descent algorithm based on the Riemannian gradient of the objective function $g(\X)$. Thus, our work has the potential to lead to new insight and extensions of these classical techniques. Part of the opportunity for novelty in this respect comes from the fact that our work considers a weighted objective function $-\frac 1 2 \text{tr}(\mathbf{X}^\top \mathbf{M} \mathbf{X}\mathbf{N})$ rather than the classical $-\frac 1 2 \text{tr}(\mathbf{X}^\top \mathbf{M} \mathbf{X})$. Part of the novelty comers from the fact that we are not limited to any one choice of algorithm for minimizing this objective function.}

\MBWbb{Second, our work also yields insight into the problem of estimating the eigenvectors of a matrix $\M$ from a randomized estimate $\Y$ of that matrix. Indeed, classical matrix perturbation theory with Davis-Kahan $\sin \Theta$ theorem~\cite[Theorem 1]{yu2015useful} can provide an upper bound\footnote{According to the Davis-Kahan $\sin \Theta$ theorem, if $\|\Y-\M\| \leq (1-1/\sqrt{2})(\lambda_r - \lambda_{r+1})$, the distance (with optimal rotation) between the two local minima can also be bounded with $\frac{2\|\Y-\M\| }{\lambda_r - \lambda_{r+1}}$~\cite{chen2020spectral}.} for the principal angle between the eigenvectors of a matrix and its perturbation. In that work, however, the eigenvectors must be arranged in a consecutive order. Our work gives more general insight into the connections between the critical points between the landscapes of the empirical and population risks. Specifically, Theorem~\ref{thm:main0} allows the critical point $\mathbf{X}_{\Omega}$ to contain any $r$ eigenvectors of $\M$ arranged in any order. (The critical point will be a strict saddle if not a global minimizer.) Corollary~\ref{coro:main} then gives an association between the critical points of the empirical and population risks. This is more general than what is provided by the Davis-Kahan $\sin \Theta$ theorem.}

\section{Applications}
\label{sec:appl}


As mentioned previously, the eigendecomposition problem~\eqref{eqn:max_eigen} is widely used in  spectral methods and can also be used to provide a ``warm start'', i.e., a good initialization to other more sophisticated algorithms, for solving various non-convex matrix optimization problems. \MBWbb{For example, the Wirtinger flow algorithm proposed in~\cite{candes2015phase} uses a spectral method to find a suitable initialization. Specifically, that paper uses the eigenvector corresponding to the largest eigenvalue of a matrix constructed from the observed measurements as the initialization for the proposed Wirtinger flow algorithm.}

In this section, we demonstrate that the two assumptions~\eqref{eqn:as_grad} and~\eqref{eqn:as_hess} can hold in a variety of applications, including matrix sensing, matrix completion, phase retrieval, and quadratic sensing. This makes it possible for us to characterize the empirical landscape of the eigendecomposition problem, where only some random measurements of $\M$ are available, via the landscape of the corresponding population risk.
In particular, one can roughly identify the positions of the local minima of the empirical risk from those of the population risk, which may result in a better understanding of the reason why the spectral initialization approaches are so powerful in providing a good initialization for these applications.
That is, instead of directly proving that the spectral initialization obtained from minimizing the empirical risk falls into the basin of attraction\footnote{Note that many iterative algorithms (e.g., gradient descent) are guaranteed to converge to a local minimum when the initialization is within the basin of attraction~\cite{candes2015phase,chi2019nonconvex}. } as in some existing literature~\cite{candes2015phase,chi2019nonconvex}, an alternative way is to show that the spectral initialization obtained from minimizing the population risk falls into the basin of attraction. Then, together with Corollary~\ref{coro:main}, one can still show that the spectral initialization obtained from minimizing the empirical risk falls into the basin of attraction.

\subsection{Matrix Sensing}
\label{sec:ms}

Consider a symmetric matrix $\M\in\RRR^{N\times N}$ with rank $r$. In the matrix sensing problem, one is given random measurements $\y\in\RRR^m$ with the $i$-th entry being
$y_i = \langle \A_i,\M \rangle, 1\leq i \leq m.$	
Here, $\{\A_i\in\RRR^{N\times N}\}_{i=1}^m$ is a set of Gaussian random sensing matrices with entries from the distribution $\NN(0,1)$. To provide a ``warm start'' with the spectral method, one can construct a surrogate matrix of $\M$ as
$\Y = \frac 1 m \sum_{i=1}^m y_i \A_i =  \frac 1 m \sum_{i=1}^m \langle \A_i,\M \rangle \A_i,$	
and then perform an eigendecomposition on $\Y$. That is, the goal is to minimize the empirical risk $f(\X)$ in~\eqref{eqn:def_er}
on the Stiefel manifold $\St(N,r)$.
Note that the expectation of $\Y$ with respect to $\A_i$ is $\EEE \Y = \M$, which implies that
$g(\X) = -\frac 1 2 \tr(\X^\top \M \X \N	)	 = \EEE f(\X),$
namely, $g(\X)$ in~\eqref{eqn:def_pr} is indeed the population risk of the empirical risk in~\eqref{eqn:def_er}.

Define a sensing operator $\AA:\RRR^{N\times N} \rightarrow \RRR^m$ with the $i$-th entry of $\AA(\X)$ being
$[\AA(\X)]_i= \frac{1}{\sqrt{m}} 	\langle \A_i,\X \rangle, ~ 1\leq i \leq m.$
As is shown in~\cite{recht2010guaranteed,candes2011tight}, the above operator $\AA$ satisfies the following Restricted Isometry Property (RIP) with probability at least $1-e^{-cm}$ if the entries of the Gaussian random sensing matrices $\A_i$ follow $\NN(0,1)$ and the number of measurements $m\geq C\delta_r^{-2} rN\log(N)$.
\begin{Definition}
(RIP~\cite{candes2008restricted}) An operator $\AA:\RRR^{N\times N} \rightarrow \RRR^m$ is said to satisfy the $r$-RIP with restricted isometry constant $\delta_r$ if
$(1-\delta_r)\|\X\|_F^2 \leq \|\AA(\X)\|_2^2 \leq (1+\delta_r)\|\X\|_F^2 $	
holds for any matrix $\X\in\RRR^{N\times N}$ with rank at most $r$.
	
\end{Definition}

Similar to the derivations in Section~\ref{sec:main}, we can write the Riemannian gradient and Hessian of $f(\X)$ as 
\begin{align*}
&\text{grad}~f(\X) = (\X\X^\top - \I_N)\Y\X \N	-  \frac 1 2 \X[\X^\top \Y \X, \N], \\
&\text{hess}~f(\X)[\U,\U] =  \langle \X^\top \Y \X, \U^\top \U \N\rangle - \langle \Y,\U\N\U^\top \rangle,
\end{align*}
where $\U\in\RRR^{N\times r}$ is a matrix that belongs to the tangent space of $\St(N,r)$ at $\X$.
Then, we have
\begin{align*}
&\sup_{\X\in \St(N,r)} \|\text{grad}~f(\X) - \text{grad}~g(\X)\|_F \\
=&\sup_{\X\in \St(N,r)} \left\|(\X\X^\top - \I_N)(\Y-\M)\X\N \right. \\
&~~ \left.- \frac 1 2 \X\X^\top (\Y-\M)\X \N + \frac 1 2 \X\N\X^\top(\Y-\M)\X \right\|_F\\
\leq& \frac 3 2 r^{\frac 3 2} \|\Y-\M\|.
\end{align*}

It follows from~\cite[Lemma 8]{chi2019nonconvex} that
$\|\Y-\M\| \leq \delta_{2r}	\sqrt{r}\|\M\|$
holds for the rank-$r$ matrix $\M$ if the sensing operator $\AA$ satisfies 2$r$-RIP with restricted isometry constant $\delta_{2r} < 1$.
Therefore, condition~\eqref{eqn:as_grad} will hold
if the RIP constant $\delta_{2r}\leq \frac{\epsilon}{3r^2\|\M\|}$.
To verify condition~\eqref{eqn:as_hess}, it suffices to show that
$|\text{hess}~f(\X)[\U,\U]-\text{hess}~g(\X)[\U,\U]|\leq \frac{\eta}{2}	$
holds for any $\U$ (with $\|\U\|_F=1$) belonging the tangent space of $\St(N,r)$ at $\X$, i.e., any $\U\in \TT_{\X}\St(N,r) = \{\U \in \RRR^{N\times r}: \X^\top \U + \U^\top \X = \zero\} = \{\X \S+ \X_\perp \K: \S^\top = -\S, \K \in \RRR^{(N-r) \times r} \}$. With Lemma 8 in~\cite{chi2019nonconvex}, we again obtain
\begin{align*}
&\sup_{\X \in \St( N,r )}      |\text{hess}~f(\X)[\U,\U]  - \text{hess}~g(\X)[\U,\U] | \\
=& \!\!	\sup_{\X\!\in \St(\!N,r\!)}\!\! |\langle \X^\top\!(\Y\!-\!\M)\X,\! \U\!^\top\!\U\N \rangle \!-\! \langle \Y\!\!-\!\M,\!\U\N\U\!^\top \!\rangle |\\
\leq & \!\!\!\sup_{\X\in \St(N,r)} \!\!\!\!\| \X\!^\top\!(\Y\!-\!\M)\X \|_F \|\U^\top\U\N  \|_F \!+\! \| (\Y\!-\!\M)\U \|_F \|\U\N\|_F\\
\leq &  2 r^{\frac 3 2}\|\Y-\M\|
\leq  \frac{\eta}{2}
\end{align*}
if the RIP constant $\delta_{2r}\leq \frac{\eta}{4r^2\|\M\|}$.


Therefore, by requiring $\delta_{2r} \leq 0.0275 r^{-2} \|\M\|^{-1} d_{\min}$ and setting $\epsilon = 3\delta_{2r} r^2 \|\M\|$ and $\eta = 0.11 d_{\min}$, we can conclude that the two conditions in~\eqref{eqn:as_grad} and~\eqref{eqn:as_hess} hold with probability at least $1-e^{-cm}$ as long as the number of measurements $m\geq Cd_{\min}^{-2} r^5 \|\M\|^2N\log(N)$. Moreover, it follows from Corollary~\ref{coro:main} that the distance between the empirical local minimum and population local minimum is on the order of $d_{\min}^{-1} \delta_{2r} r^2 \|\M\|$.

\subsection{Matrix Completion}
\label{sec:mc}

In this scenario, one is given measurements $\PP_{\Omega}(\M)$, where $\PP_\Omega:\RRR^{N\times N}\rightarrow \RRR^{N\times N}$ is a projection operator defined with
$(\PP_{\Omega}(\M))_{i,j} = M_{i,j}$ if $(i,j)\in\Omega$ and $0$ otherwise.
Assume that each entry of a positive semidefinite (PSD) matrix $\M\in \RRR^{N\times N}$ is observed independently with probability $0<p\leq 1$, namely, $(i,j)\in\Omega$ independently with probability $p$. Suppose that $\M$ is a low-rank matrix with rank $r$. Define $\Ut \triangleq [\sqrt{\lambda_1}\u_1~\cdots~\sqrt{\lambda_r}\u_r]$, where $\lambda_i$ is the $i$-th largest eigenvalue of $\M$ and $\u_i$ is the corresponding eigenvector. We also assume that $\M$ is $\mu_r$-incoherent with $\mu_r = \frac{N}{r} \max_i \sum_{j=1}^r \widetilde{U}_{i,j}^2$ being well-bounded~\cite{candes2009exact,keshavan2010matrix,sun2016guaranteed,chen2019model} .

To provide a ``warm start'' with the spectral method, one can construct a surrogate matrix of $\M$ as
$\Y =\frac1 p \PP_\Omega(\M) = \sum_{(i,j)\in\Omega} \langle \A_{i,j},\M \rangle \A_{i,j}$
with $\A_{i,j} = \frac{1}{\sqrt{p}} \e_i \e_j^\top$, where $\e_i$ denotes the $i$-th column of an $N\times N$ identity matrix $\I_N$. Then, performing eigendecomposition on $\Y$ corresponds to minimizing the  empirical risk $f(\X)$ in~\eqref{eqn:def_er}
on the Stiefel manifold $\St(N,r)$. Note that the expectation of $\Y$ is $\EEE \Y = \M$, which implies that
$g(\X) = -\frac 1 2 \tr(\X^\top \M \X \N	)	 = \EEE f(\X),$
which further implies that the function $g(\X)$ in~\eqref{eqn:def_pr} is also the population risk of the empirical risk in~\eqref{eqn:def_er}.

We again obtain
\begin{equation}
\begin{aligned}
&\sup_{\X\in \St(N,r)}\!\!\! \|\text{grad}~f(\X) \!-\! \text{grad}~g(\X)\|_F \!\leq \! \frac 3 2 r^{\frac 3 2} \|\Y-\M\|, \\
&\sup_{\X\in \St(N,r)}\!\! |\text{hess}~f(\X)[\U,\!\U]\!-\!\text{hess}~g(\X)[\U,\!\U]|
 \!\leq\! 2 r^{\frac 3 2}\!\|\Y\!-\!\M\|.
\label{eqn:bound_2as}
\end{aligned}
\end{equation}
It follows from~\cite[Theorem 6.3]{candes2009exact} that
\begin{align*}
\|\Y \!-\! \M\|  \!=\!  p ^{-1} \| (\PP_\Omega  \!-\!  p\mathcal{I})(\M) \|
 \!\leq\!  C  (Np) ^{-\frac 1 2}\! \sqrt{\log(N)}r\lambda_1 \mu_r
\end{align*}
holds with probability at least $1-N^{-c_1}$ if $p\geq c_2 \frac{\log(N)}{N}$. The last inequality follows from the fact that $	\|\M\|_{\infty} \leq 	\frac{r}{N}\lambda_1 \mu_r$~\cite{chen2019model}.


Therefore, by setting $\epsilon = Cr^{\frac 5 2} \lambda_1 \mu_r (Np) ^{-\frac 1 2} \sqrt{\log(N)} \leq \eta = 0.11 d_{\min}$,  we can conclude that the two conditions in~\eqref{eqn:as_grad} and~\eqref{eqn:as_hess} hold with probability at least $1-N^{-c_1}$ as long as $p\geq C \max\{N^{-1} \log(N), r^5 \lambda_1^2 \mu_r^2 d_{\min}^{-2} N^{-1} \log(N)  \}$. Moreover, it follows from Corollary~\ref{coro:main} that the distance between the empirical local minimum and population local minimum is on the order of $d_{\min}^{-1} r^{\frac 5 2} \lambda_1 \mu_r (Np) ^{-\frac 1 2} \sqrt{\log(N)}$.

\subsection{Phase Retrieval}
\label{sec:pr}

In phase retrieval, one is given the  measurements of a vector $\xs\in\RRR^N$:
$y_i = (\a_i^\top \xs)^2, 1\leq i \leq m,$
where $\{\a_i\}$ are sensing vectors with entries following $\NN(0,\frac{1}{\sqrt{2}})$.
Denote $\M = \xs{\xs}^\top$. Note that the above measurements can be rewritten as
$y_i = {\a_i}^\top \M \a_i= \langle \a_i\a_i^\top,\M \rangle,  1\leq i \leq m,$
which can be viewed as a special case of the matrix sensing problem with the sensing matrix $\A_i = \a_i\a_i^\top$.

To provide a ``warm start'' with the spectral method, one can construct a surrogate matrix of $\M$ as
$\Y = \frac 1 m \sum_{i=1}^m y_i \a_i \a_i^\top,$	
and then perform eigendecomposition on $\Y$, minimizing the following empirical risk
$f(\x) = -\frac 1 2 \x^\top \Y \x$	
on the unit sphere $\SSS^{N-1}$. Note that the expectation of $\Y$ with respect to $\a_i$ is $\EEE \Y = \M + \frac 1 2 \|\xs\|_2^2  \I_N$, which implies that
$g(\x) = -\frac 1 2 \x^\top \M \x	 -\frac 1 4  \|\xs\|_2^2 \|\x\|_2^2 = \EEE f(\x),$
namely, the above function $g(\x)$ is indeed the population risk of the empirical risk $f(\x)$.


With some elementary calculations, we write the Riemannian gradient and Hessian of the empirical risk $f(\x)$ as 
\begin{align*}
&\text{grad}~f(\x) = (\x\x^\top\! \!-\! \I_N)(\Y\x\!+\!\frac 1 2 \|\xs\|_2^2 \x), \\
&\text{hess}~f(\x) = (\I_N\!-\!\x\x^\top\!) (\x^\top\Y\x\I_N\!-\!\Y) (\I_N\!-\!\x\x^\top\!).
\end{align*}
Then, similar to the matrix sensing case, we have
\begin{align*}
&\sup_{\x\in \SSS^{N-1}} \|\text{grad}~f(\x) - \text{grad}~g(\x)\|_2 \\
=&\sup_{\|\x\|_2 = 1} \|(\x\x^\top - \I_N)(\Y-\M-\frac 1 2 \|\xs\|_2^2 \I_N)\x\|_2\\
\leq & 	\sup_{\|\x\|_2 = 1} \|\x\x^\top - \I_N\| \|\Y-\EEE \Y\|\|\x\|_2
= \|\Y-\EEE\Y\|,
\end{align*}
and
\begin{align*}
&\sup_{\x\in \SSS^{N-1}} |\u^\top \text{hess}~f(\x)\u-\u^\top \text{hess}~g(\x)\u|\\
= 	&\sup_{\|\x\|_2=1} |\x^\top (\Y-\EEE\Y)\x - \u^\top(\Y-\EEE\Y)\u|\\
\leq & \sup_{\|\x\|_2=1} |\x^\top (\Y-\EEE\Y)\x| + |\u^\top(\Y-\EEE\Y)\u|\\
\leq &  2\|\Y-\EEE\Y\|.
\end{align*}

For any small constant $\delta \in(0,1]$, it follows from~\cite[Lemma 13]{li2019nonconvex} that
$
\|\Y-\EEE\Y\| \leq \delta \|\xs\|_2^2	
$
with probability at least $1-c_1N^{-c_2}$ provided that $m \geq C \delta^{-2} N\log(N)$ for some constant $C$.
Therefore, by setting $\epsilon = C\delta \|\xs\|_2^2 \leq \eta = 0.11 d_{\min}\SLbb{=0.11 \|\x^\star\|_2^2}$ \MBWbb{(since $d_{\min} = \|\x^\star\|_2^2$)}, we can conclude that the two conditions in~\eqref{eqn:as_grad} and~\eqref{eqn:as_hess} hold with probability at least $1-c_1N^{-c_2}$ as long as the number of measurements satisfies $m\geq C d_{\min}^{-2} \|\xs\|_2^4 N \log(N) \SLbb{= C N \log(N)}$. Moreover, it follows from Corollary~\ref{coro:main} that the distance between the empirical local minimum and population local minimum is on the order of $d_{\min}^{-1} \delta \|\xs\|_2^2 \SLbb{=\delta}$.

\subsection{Quadratic Sensing}
\label{sec:qs}

The quadratic sensing problem can be viewed as a rank-$r$ extension of the phase retrieval problem. In particular, consider a PSD matrix $\M=\Xs{\Xs}^\top$ with $\Xs\in\RRR^{N\times r}$, and suppose one is given the measurements
$y_i = {\a_i}^\top \M \a_i= \langle \a_i\a_i^\top,\M \rangle,  1\leq i \leq m,	$
where $\{\a_i\}$ are sensing vectors with entries following $\NN(0,\frac{1}{\sqrt{2}})$.
Again, to provide a ``warm start'' with the spectral method, one can construct a surrogate matrix of $\M$ as
$\Y = \frac 1 m \sum_{i=1}^m y_i \a_i \a_i^\top,$	
and then perform eigendecomposition on $\Y$, i.e., minimizing the empirical risk $f(\X)$ in~\eqref{eqn:def_er}
on the Stiefel manifold $\St(N,r)$. Note that the expectation of $\Y$ with respect to $\a_i$ is $\EEE \Y = \M + \frac 1 2 \|\Xs\|_F^2 \I_N$, which implies that
$g(\X) = -\frac 1 2 \tr(\X^\top \M \X \N	)	 -\frac 1 8 r(r+1) \|\Xs\|_F^2  = \EEE f(\X),$
i.e., the above function $g(\X)$ is again the population risk of the empirical risk $f(\X)$ in~\eqref{eqn:def_er}.


Similar to the other three applications, we can again bound the left hand side of the two assumptions~\eqref{eqn:as_grad} and~\eqref{eqn:as_hess} with $ r^{\frac 3 2} \|\Y-\EEE\Y\|$. Then,
for any small constant $\delta \in(0,1]$, it follows from~\cite[Corollary 5.2]{sanghavi2017local} that
$\|\Y-\EEE\Y\| \leq \delta \|\Xs\|_F^2$	
with probability at least $1-c_1e^{-c_2rN}-c_3m^{-2}$ when provided $m \geq C \delta^{-2} rN\log^2(N)$ for some constant $C$.

Therefore, by setting $\epsilon = C\delta r^{\frac 3 2} \|\Xs\|_F^2 \leq \eta = 0.11 d_{\min}$, we can conclude that the two conditions in~\eqref{eqn:as_grad} and~\eqref{eqn:as_hess} hold with probability at least $1-c_1e^{-c_2rN}-c_3m^{-2}$ as long as the number of measurements satisfies $m \geq C  d_{\min}^{-2} r^4 \|\Xs\|_F^4 N \log^2(N)$. Moreover, it follows from Corollary~\ref{coro:main} that the distance between the empirical local minimum and population local minimum is on the order of $d_{\min}^{-1} \delta r^{\frac 3 2 } \|\Xs\|_F^2$.


%
%
%

\section{Simulation Results}
\label{sec:simu}

In this section, we conduct experiments to further support our theory. In particular, we illustrate our theory using the four fundamental low-rank matrix optimization problems that are detailed in Section~\ref{sec:appl}.

\begin{figure}[t]
\begin{minipage}{0.22\linewidth}
\centering
\includegraphics[width=0.8in]{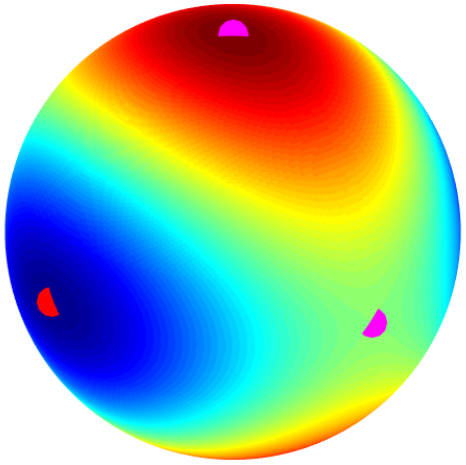}
\end{minipage}
\hfill
\begin{minipage}{0.22\linewidth}
\centering
\includegraphics[width=0.8in]{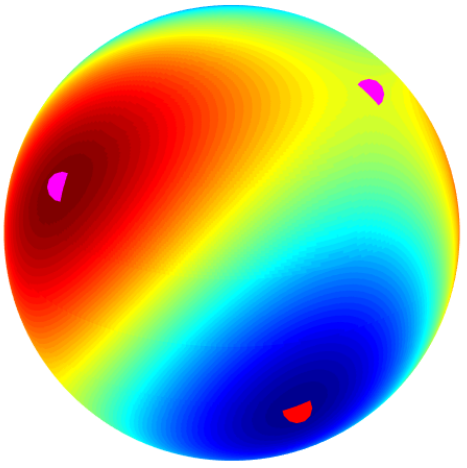}
\end{minipage}
\hfill
\begin{minipage}{0.22\linewidth}
\centering
\includegraphics[width=0.8in]{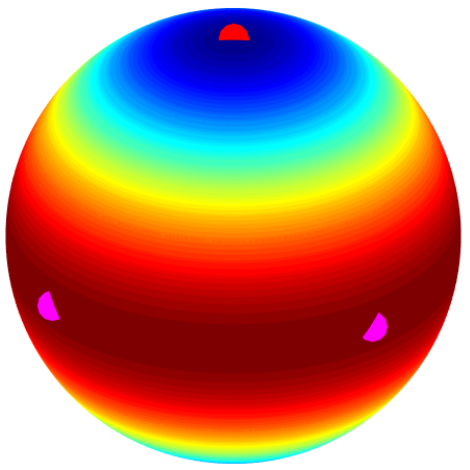}
\end{minipage}
\hfill
\begin{minipage}{0.22\linewidth}
\centering
\includegraphics[width=0.8in]{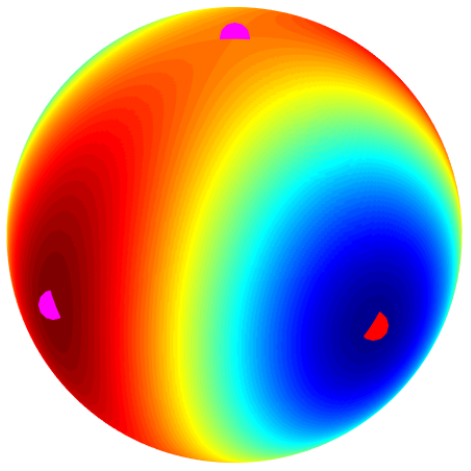}
\end{minipage}\\
\begin{minipage}{0.22\linewidth}
\centering
\includegraphics[width=0.77in]{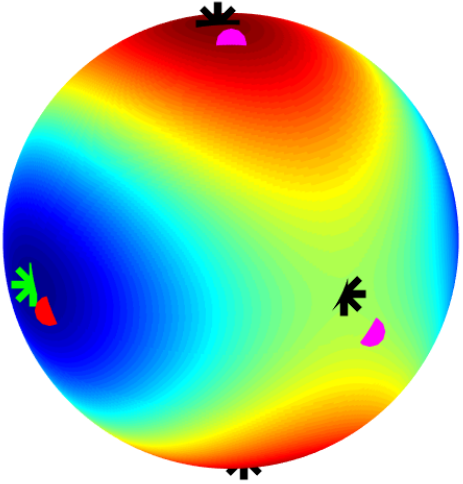}
\centerline{\scriptsize{(a) matrix sensing}}
\end{minipage}
\hfill
\begin{minipage}{0.22\linewidth}
\centering
\includegraphics[width=0.8in]{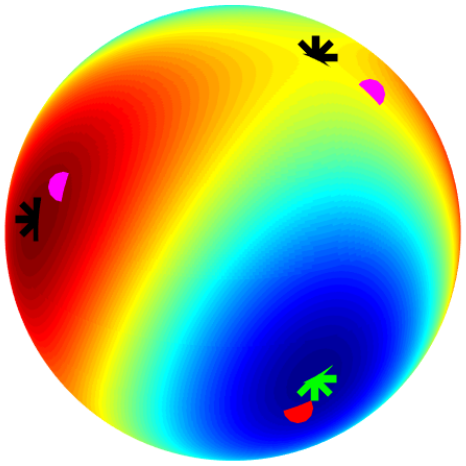}
\centerline{\scriptsize{(b) matrix completion}}
\end{minipage}
\hfill
\begin{minipage}{0.22\linewidth}
\centering
\includegraphics[width=0.8in]{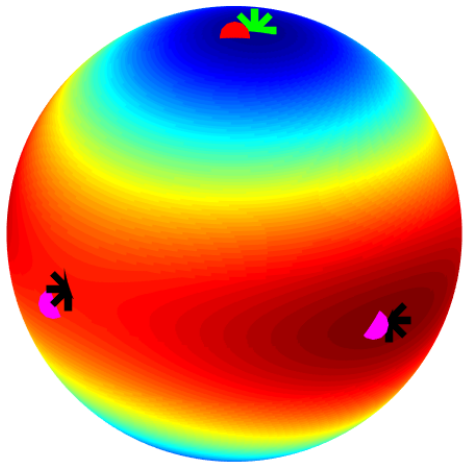}
\centerline{\scriptsize{(c) phase retrieval}}
\end{minipage}
\hfill
\begin{minipage}{0.22\linewidth}
\centering
\includegraphics[width=0.8in]{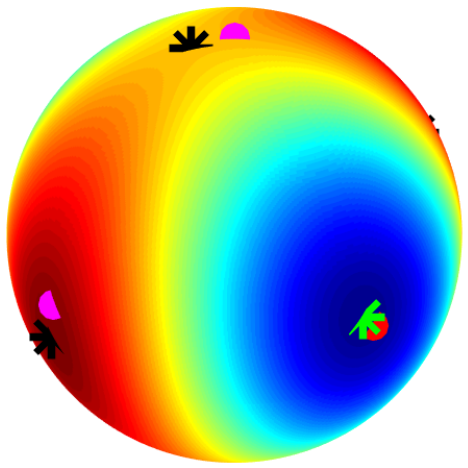}
\centerline{\scriptsize{(d) quadratic sensing}}
\end{minipage}
\caption{The landscape of population risk and empirical risk of the eigendecomposition problem used in spectral methods for a variety of low-rank matrix optimization problems. The top four figures represent the landscape of the population risk, and the bottom four figures represent a realization of the empirical risk. We use red and magenta dots to indicate global minima and (strict) saddle points of the population risks, respectively, and we use green and black stars to indicate global minima and (strict) saddle points of the empirical risks, respectively. The blue color indicates smaller objective values and the red color indicates larger values.}
\label{fig:landscape}
\end{figure}

In matrix sensing, for visualization purposes, we use a small true data matrix $\M=\text{diag}([2,1,0])$ (i.e., $N=3$) and obtain $m=100$ samples with the measurement model introduced in Section~\ref{sec:ms}. In matrix completion (see Section~\ref{sec:mc}), we set the true data matrix as $\M=\U\U^\top$, where $\U \in\RRR^{3\times 2}$ is generated as a Gaussian random matrix with normalized columns. We set the sampling probability as $p=0.8$. In phase retrieval, we set $\xs = [0;0;1]$ and then construct $m=200$ samples with the measurement model introduced in Section~\ref{sec:pr}. Finally, in quadratic sensing, we construct the data matrix as $\M=\Xs{\Xs}^\top$ with $\Xs = [0~0;0~2;1~0]\in\RRR^{3\times 2}$. Then, we construct $m=100$ samples based on the measurement model introduced in Section~\ref{sec:qs}.
Surrogate matrices $\Y$ for these applications are then constructed with the given random samples. We present the landscapes of the population risk and empirical risk of the eigendecomposition problem used in the spectral method for these low-rank matrix optimization problems in Figure~\ref{fig:landscape}. It can be seen that there does exist a direct correspondence between the local (global) minima of the empirical risk and population risk in each of these applications when the number of measurements is sufficiently high.\footnote{Note that we need a relatively large number of measurements because of the small size of the data matrix.
}

\begin{figure}[t]
\begin{minipage}{0.48\linewidth}
\centering
\includegraphics[width=1.67in]{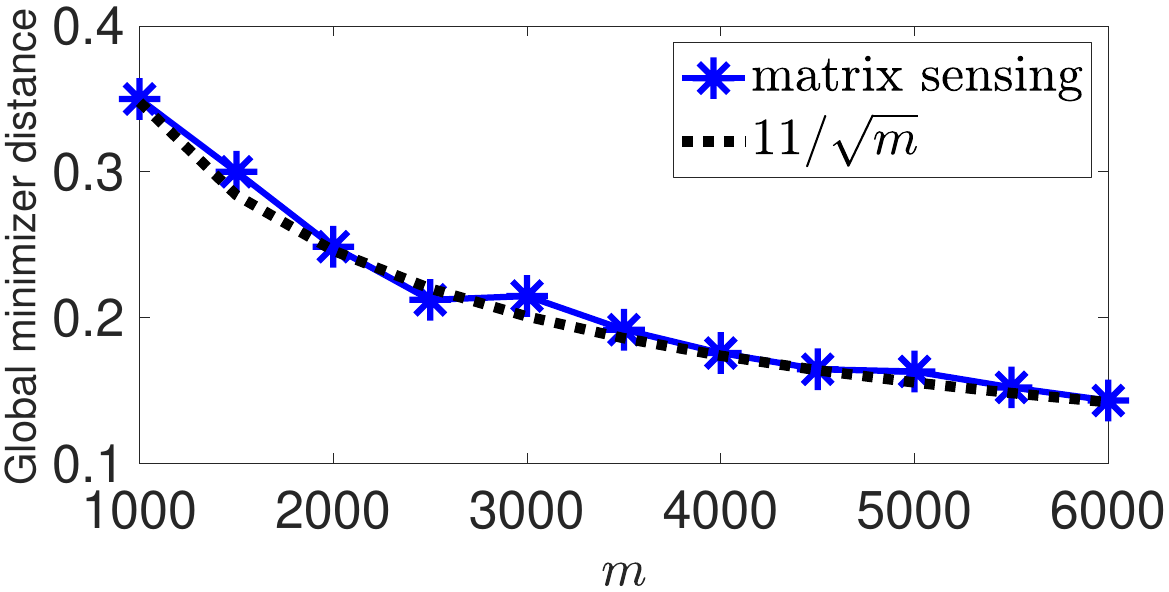}
\centerline{\footnotesize{(a) matrix sensing}}
\end{minipage}
\hfill
\begin{minipage}{0.48\linewidth}
\centering
\includegraphics[width=1.67in]{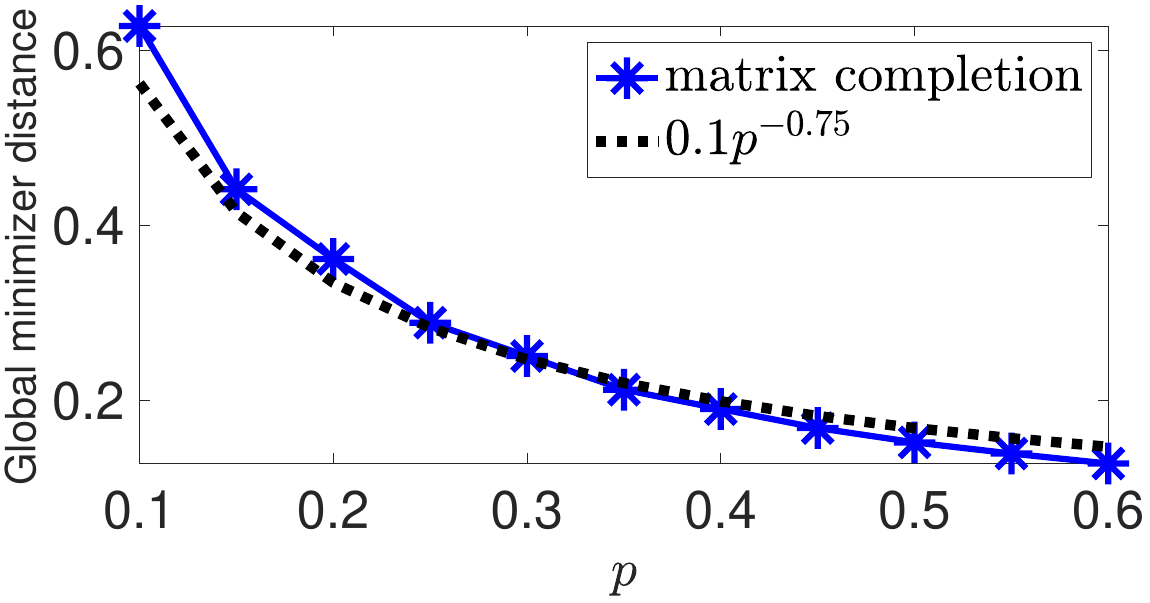}
\centerline{\footnotesize{(b) matrix completion}}
\end{minipage}
\\
\begin{minipage}{0.48\linewidth}
\centering
\includegraphics[width=1.67in]{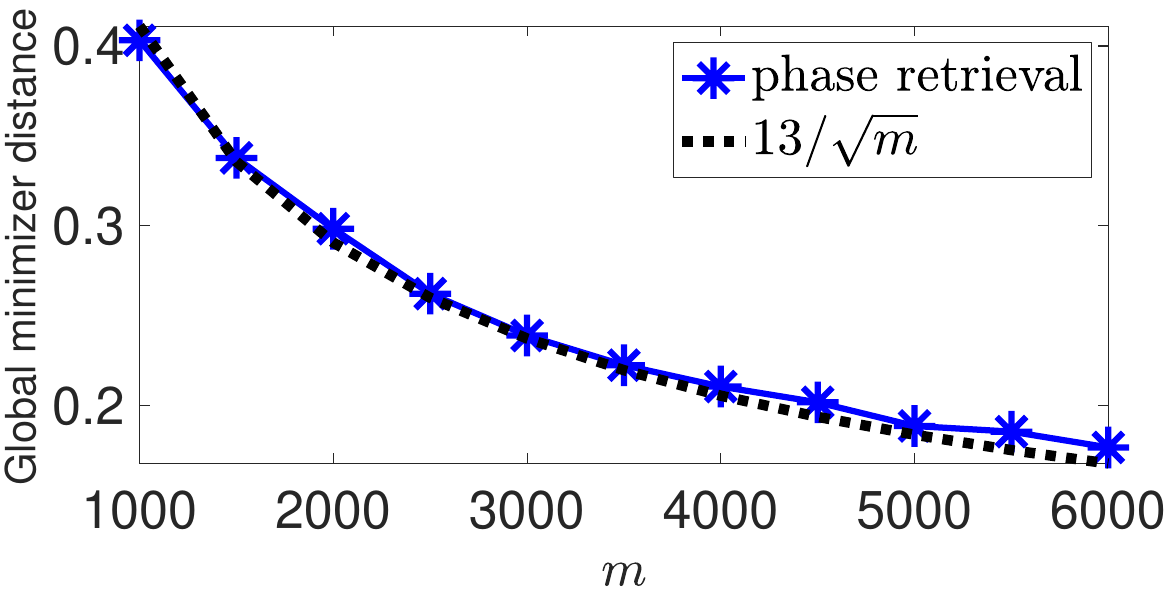}
\centerline{\footnotesize{(c) phase retrieval}}
\end{minipage}
\hfill
\begin{minipage}{0.48\linewidth}
\centering
\includegraphics[width=1.67in]{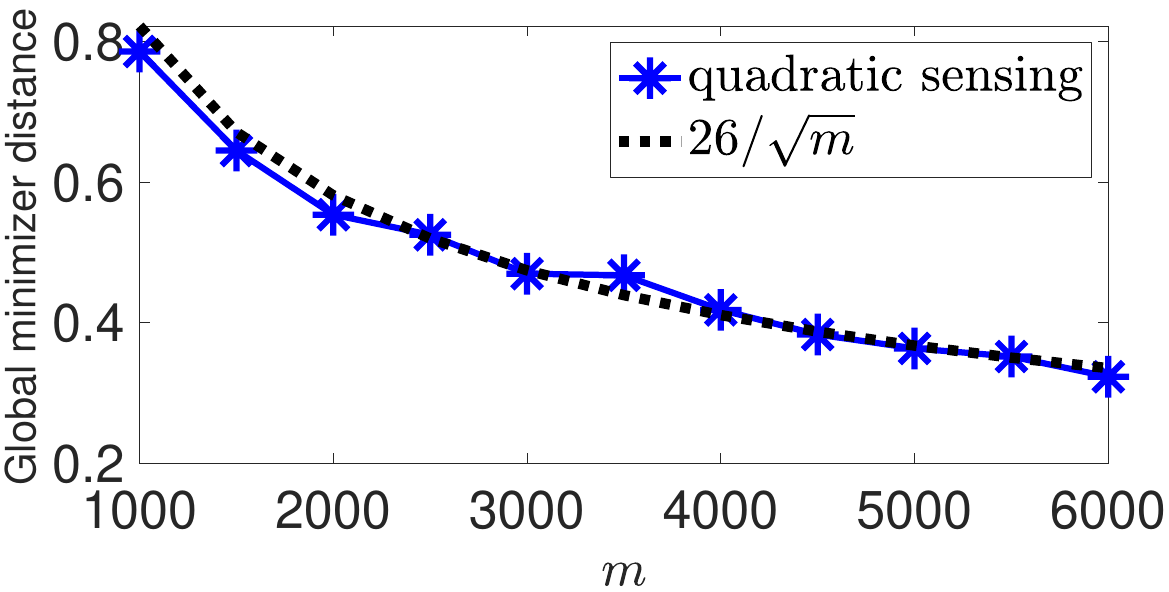}
\centerline{\footnotesize{(d) quadratic sensing}}
\end{minipage}
\caption{The distance between the population and empirical global minimizers with respect to different numbers of samples.}
\label{fig:dist}
\end{figure}

Next, we test our results in a higher-dimensional setting with $N=100$ and illustrate how the distance between the population and empirical global minimizers scales with the number of samples. In matrix sensing, matrix completion, and quadratic sensing, we set the true data matrix as $\M = \U \U^\top$, where $\U \in \RRR^{100\times 2}$ is generated as a Gaussian random matrix with normalized columns. In phase retrieval, we generate $\xs$ as a length-$N$ random vector and normalize it before we construct the random samples. We present the distance between the population and empirical global minimizers with respect to different numbers of samples in Figure~\ref{fig:dist}. Note that the results are averaged over 20 trials. It can be seen that the distance roughly scales with $p^{-3/4}$ in matrix completion and $1/\sqrt{m}$ in the other three applications, which \SLbb{is consistent}
with the analysis in Section~\ref{sec:appl}. Precisely, as is shown at the end of Section~\ref{sec:ms}, the distance between the empirical local minimum and population local minimum in matrix sensing scales with $\delta_{2r}$, and together with $m\geq C\delta_{2r}^{-2} rN\log(N)$, one can conclude that the distance scales with $1/\sqrt{m}$. A similar analysis can be conducted in phase retrieval and quadratic sensing. In matrix completion, we have shown at the end of Section~\ref{sec:mc} that the distance scales with $p^{-0.5}$, which is slightly loose when compared with the numerical observations.

\begin{figure}[t]
\begin{minipage}{0.48\linewidth}
\centering
\includegraphics[width=1.67in]{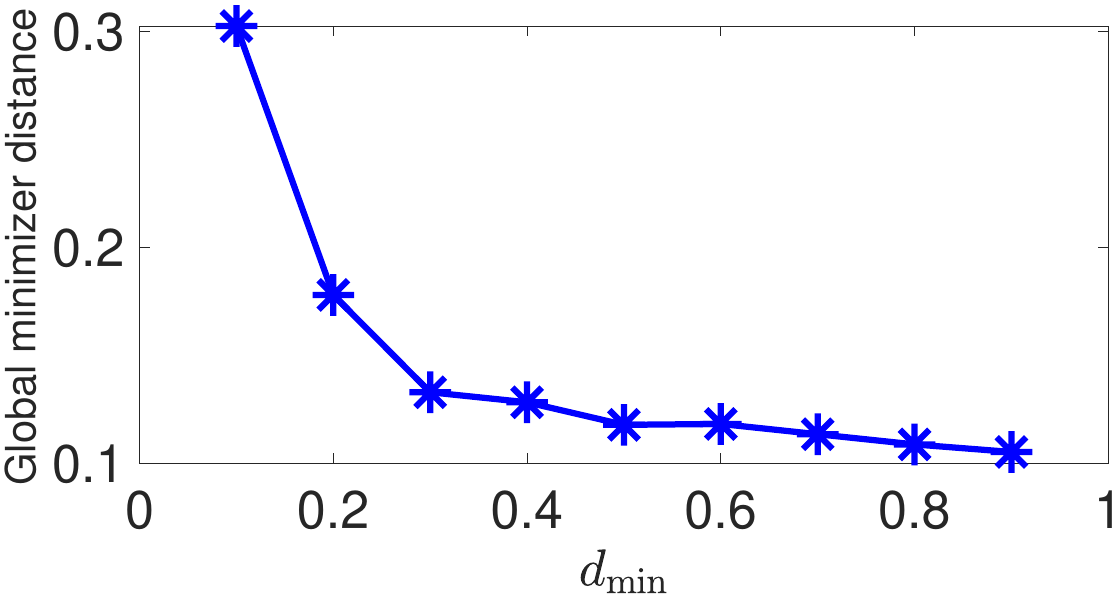}
\centerline{\footnotesize{(a) matrix sensing}}
\end{minipage}
\hfill
\begin{minipage}{0.48\linewidth}
\centering
\includegraphics[width=1.67in]{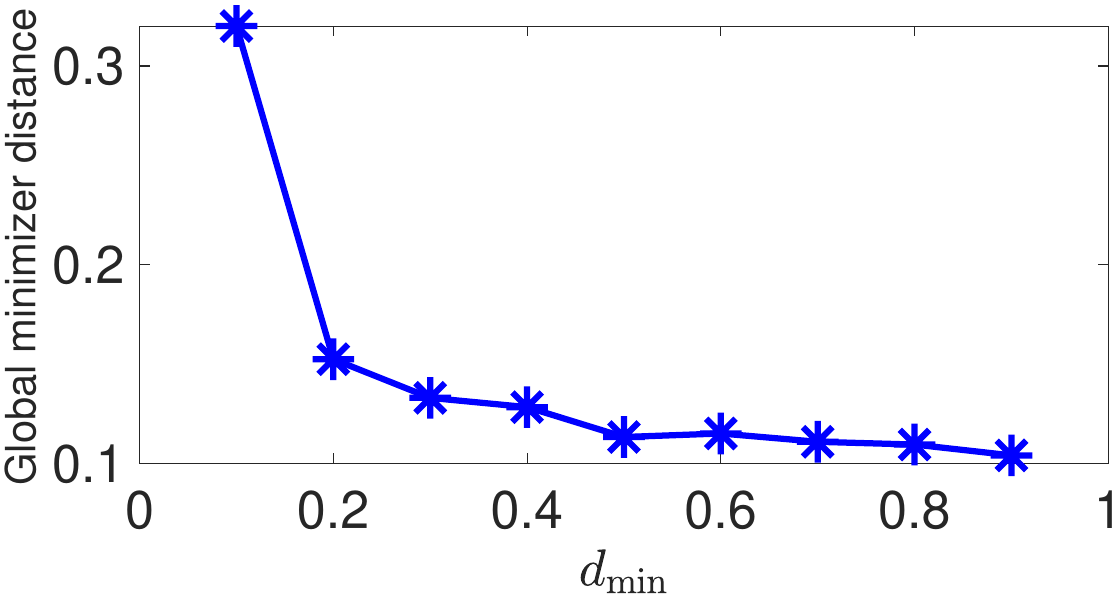}
\centerline{\footnotesize{(b) matrix completion}}
\end{minipage}
\\
\begin{minipage}{0.48\linewidth}
\centering
\includegraphics[width=1.67in]{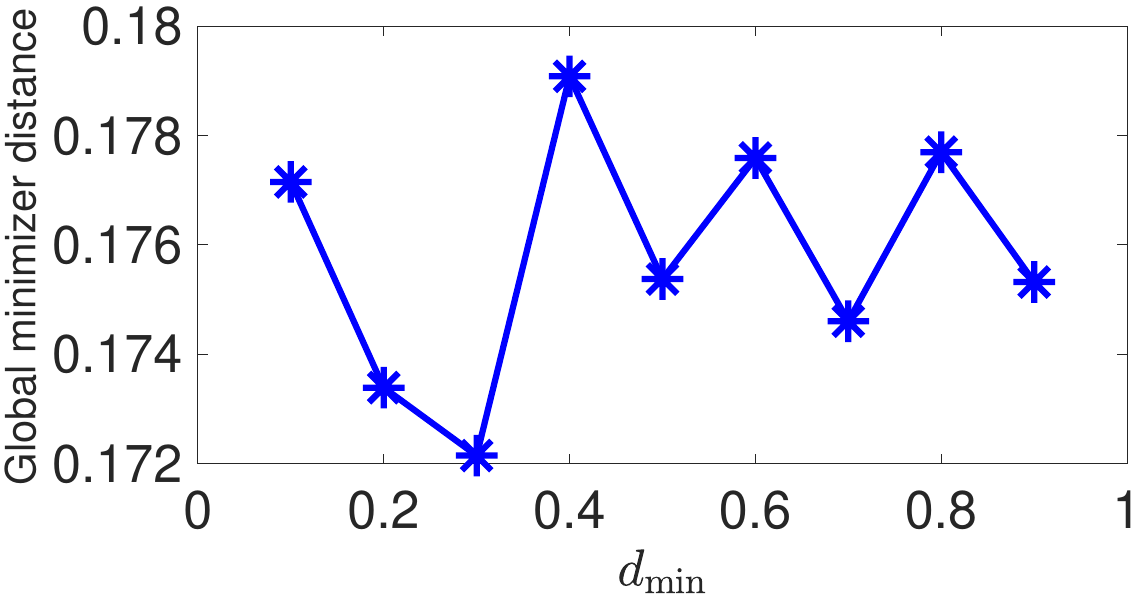}
\centerline{\footnotesize{(c) phase retrieval}}
\end{minipage}
\hfill
\begin{minipage}{0.48\linewidth}
\centering
\includegraphics[width=1.67in]{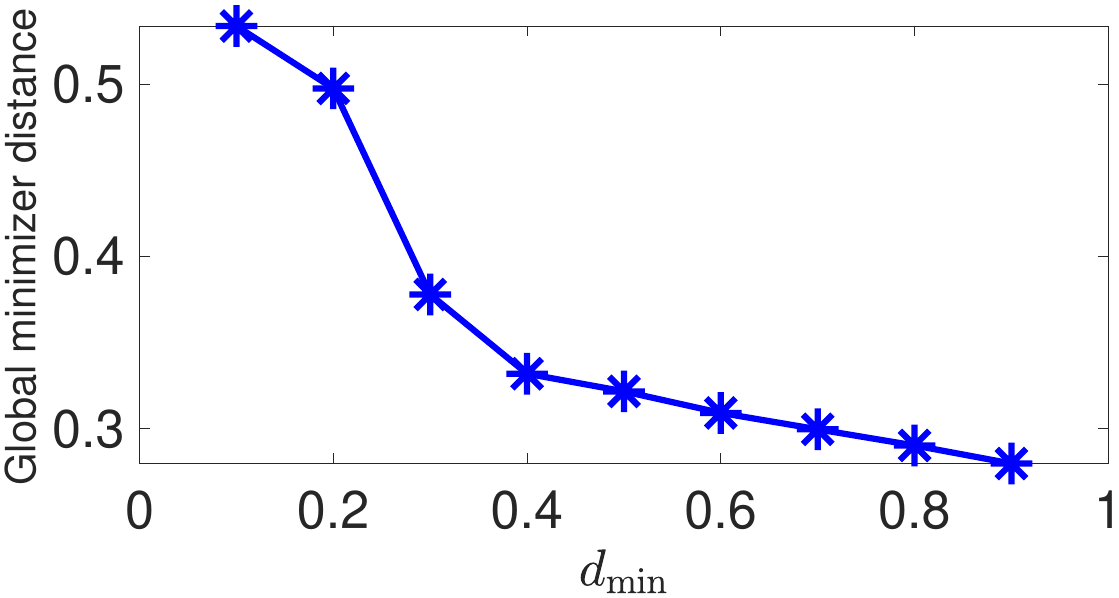}
\centerline{\footnotesize{(d) quadratic sensing}}
\end{minipage}
\caption{The distance between the population and empirical global minimizers with respect to different $d_{\min}$.}
\label{fig:dist_dmin}
\end{figure}

\SLbb{Finally, we repeat the above experiments to illustrate how the distance between the population and empirical global minimizers scales with $d_{\min}$, the minimal distance between any two of the first $r+1$ eigenvalues. We take 6000 random samples in each application. In matrix sensing, matrix completion, and quadratic sensing, we set the true data matrix as $\M = \U \Lambda \U^\top$, where $\U \in \RRR^{100\times 2}$ is generated as a Gaussian random matrix with orthonormal columns and $\Lambda \in \RRR^{2 \times 2}$ is a diagonal matrix with two diagonal entries $\lambda_1 = 2$ and $\lambda_2 = \lambda_1 - d_{\min}$.
In phase retrieval, we generate $\xs$ as a length-$N$ random vector with $\|\xs\|_2^2 = d_{\min}$. The other settings are same as the above experiments. It can be seen from Figure~\ref{fig:dist_dmin} that the distance between the population and empirical global minimizers does behave like a constant in phase retrieval and scales inversely with $d_{\min}$ in the other three applications, which is consistent with our theory.
}

\section{Conclusion}
\label{sec:conc}

In this work, we study the landscape of the eigendecomposition problem that is widely used in spectral methods. In particular, we generalize the existing analysis of the landscape of the eigendecomposition problem at the critical points to larger regions near the critical points in a special case of finding the first leading eigenvector, and we extend these results to a more general eigendecomposition problem, i.e., finding the first $r$ leading eigenvectors. Moreover, we build a connection between the landscape of the eigendecomposition problem using random measurements (empirical risk) and that of the problem using the true data matrix (population risk). With this connection, one may analyze the landscape of the empirical risk using the landscape of the population risk, which could also lead to a better understanding of why these spectral initialization approaches are so powerful.

\section*{Acknowledgement}

\SLb{MW and SL were supported by NSF grant CCF-1704204.}

\ifCLASSOPTIONcaptionsoff
  \newpage
\fi

\bibliographystyle{ieeetr}
\bibliography{SS_eigvec}

\begin{IEEEbiography}[{\includegraphics[width=1in,height=1.25in,clip,keepaspectratio]{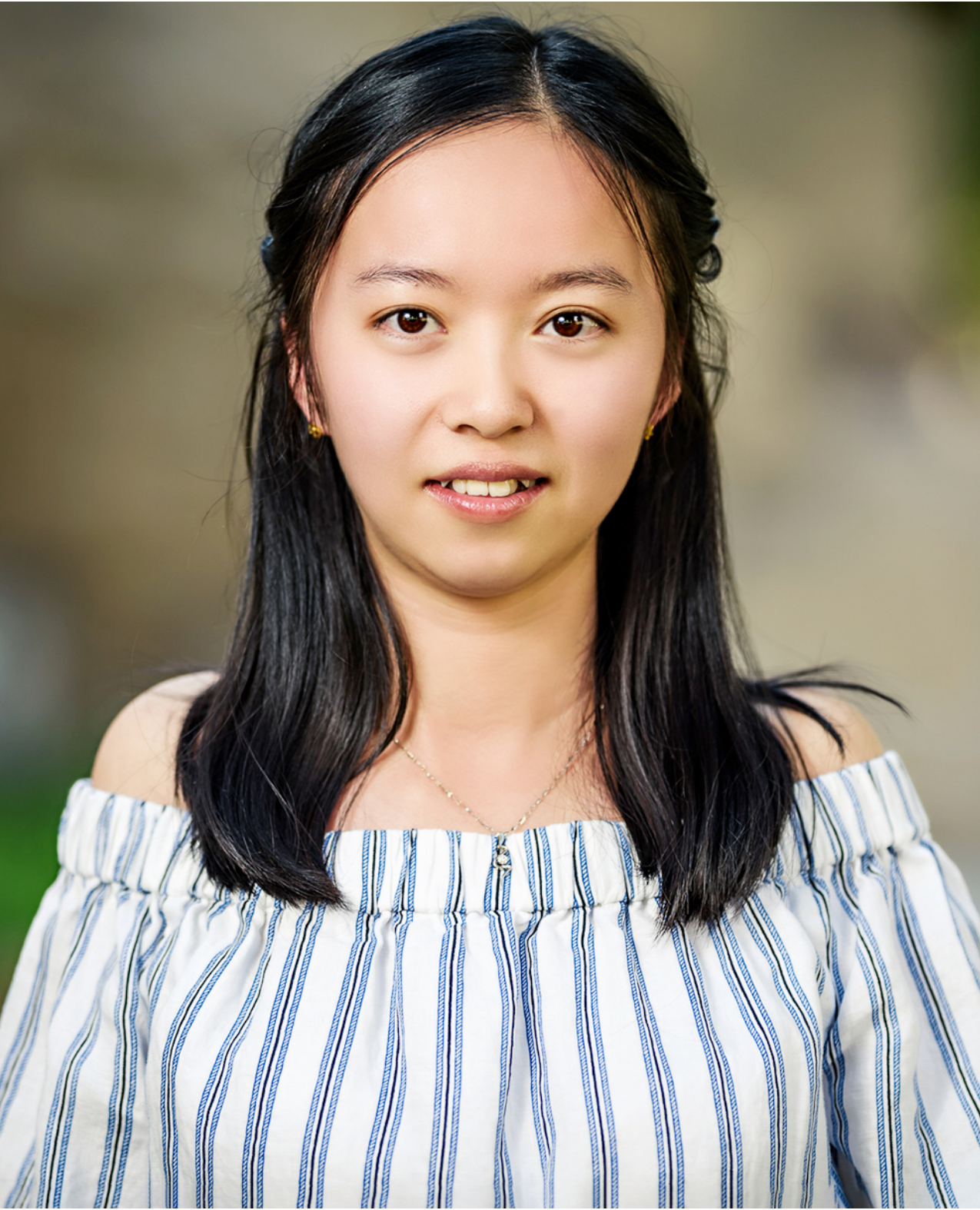}}]{Shuang Li}
received the B. Eng. degree in communications engineering from Zhejiang University of Technology, Hangzhou, China, in 2013, and the Ph.D. degree in electrical engineering from the Colorado School of Mines, Golden, CO, USA, in 2020.

She is currently a Hedrick Assistant Adjunct Professor with the Department of Mathematics, University of California, Los Angeles, CA, USA.
Her research interests include developing optimization-based techniques with optimality guarantees for fundamental problems in signal processing and machine learning.
\end{IEEEbiography}

\begin{IEEEbiography}[{\includegraphics[width=1in,height=1.25in,clip,keepaspectratio]{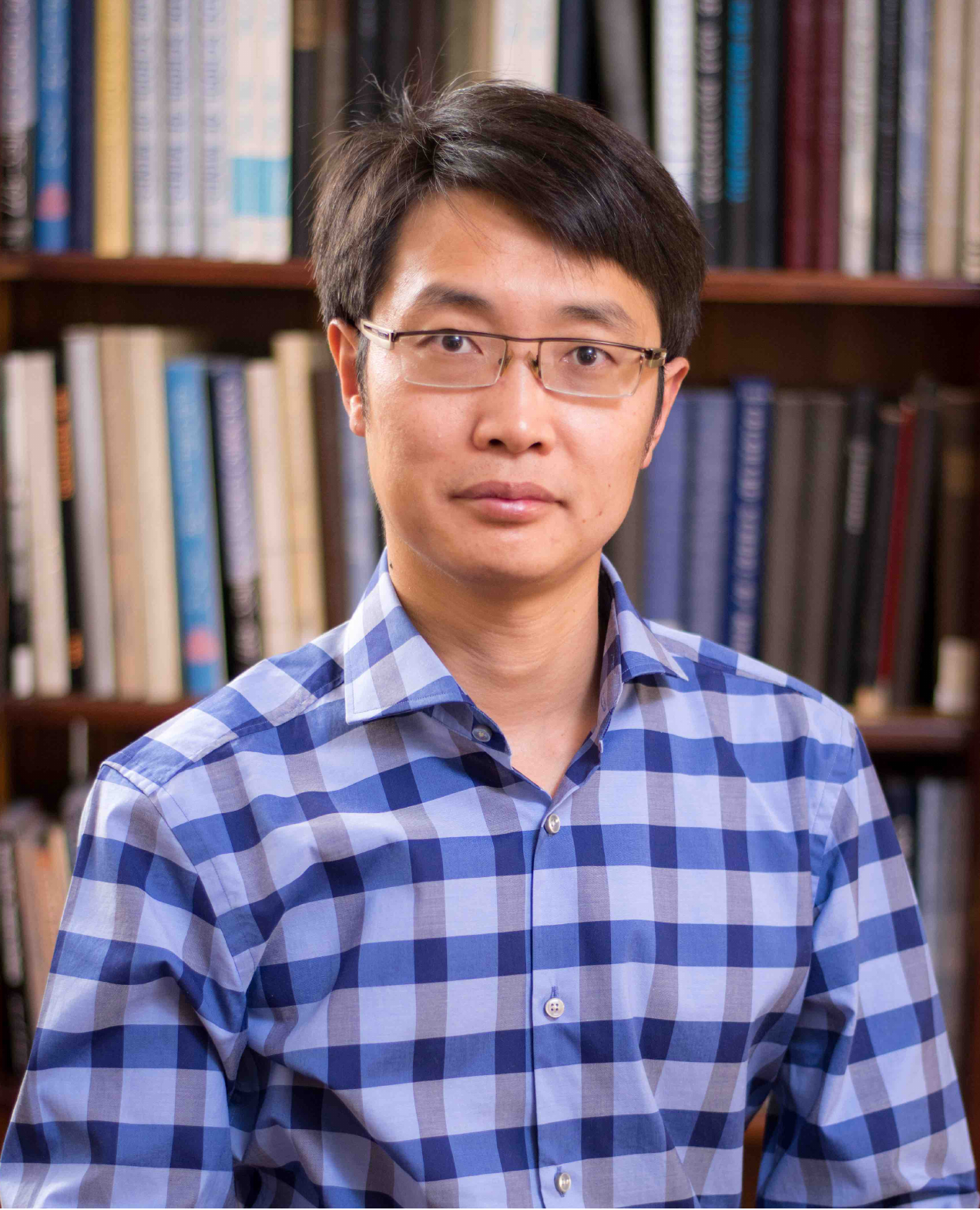}}]{Gongguo Tang}
(S'09-M'11) started his academic career as an Assistant Professor at the Colorado School of Mines in 2014, and was tenured and promoted to associate professor in 2020. He was a postdoc at the University of Wisconsin-Madison and a visiting scholar to the big data program at the Simons Institute, University of California-Berkeley from 2011 to 2013. He received his PhD in electrical engineering from Washington University in St. Louis in 2011.

He is now an Associate Professor in the Department of Electrical, Computer \& Energy Engineering at the University of Colorado-Boulder. His research revolves around modeling and optimization to extract information from data through computation. He is especially interested in the design of learning models, optimization formulations, and numerical procedures that come with theoretical performance guarantees and are scalable to large datasets, with target applications in signal processing, machine learning, and imaging.
\end{IEEEbiography}

\begin{IEEEbiography}[{\includegraphics[width=1in,height=1.25in,clip,keepaspectratio]{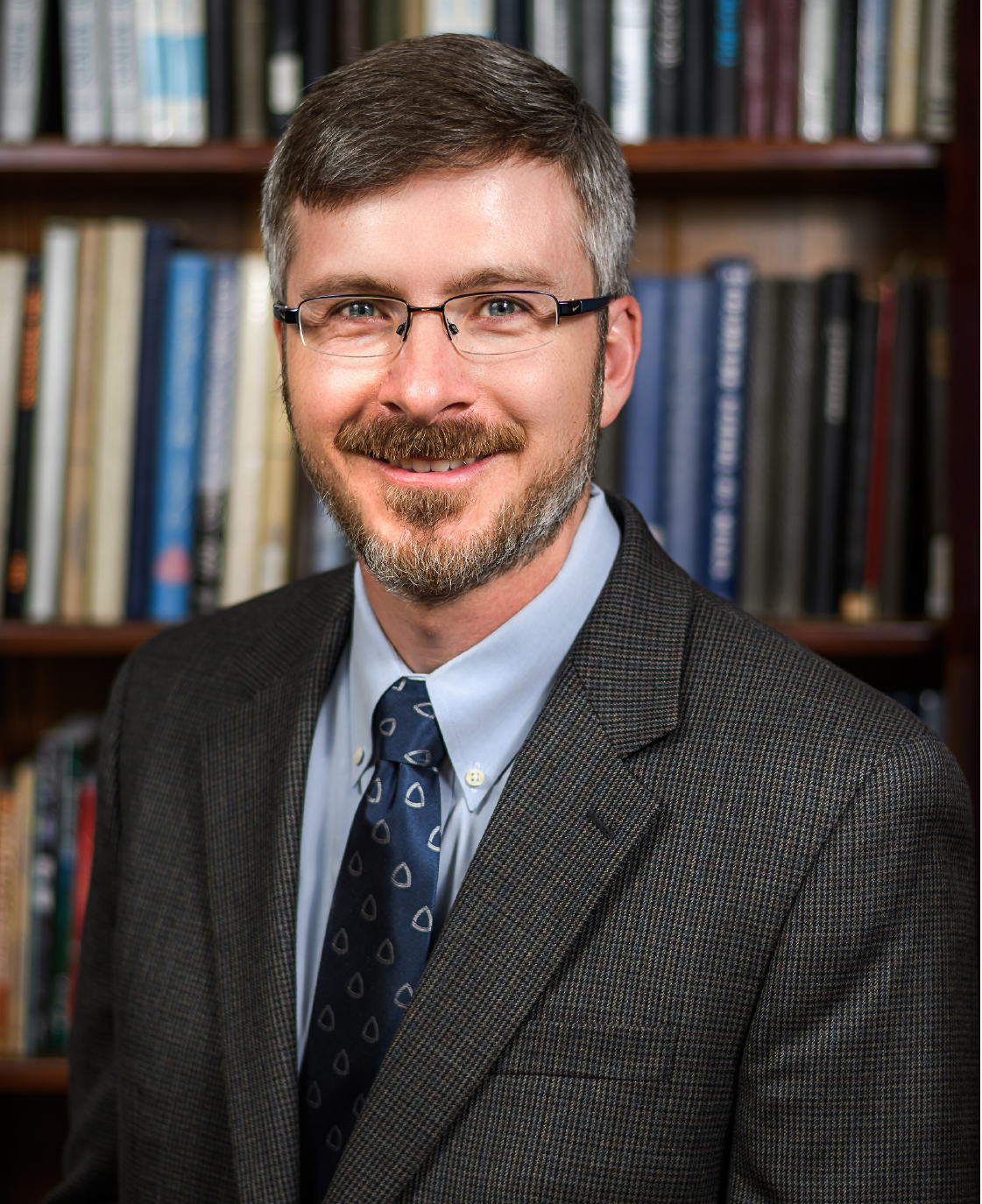}}]{Michael B. Wakin}
(S'01-M'06-SM'13-F'21) is a Professor of Electrical Engineering at the Colorado School of Mines. Dr. Wakin received a B.S. in electrical engineering and a B.A. in mathematics in 2000 (summa cum laude), an M.S. in electrical engineering in 2002, and a Ph.D. in electrical engineering in 2007, all from Rice University. He was an NSF Mathematical Sciences Postdoctoral Research Fellow at Caltech from 2006-2007, an Assistant Professor at the University of Michigan from 2007-2008, and a Ben L. Fryrear Associate Professor at Mines from 2015-2017. His research interests include signal and data processing using sparse, low-rank, and manifold-based models.

In 2007, Dr. Wakin shared the Hershel M. Rich Invention Award from Rice University for the design of a single-pixel camera based on compressive sensing. In 2008, Dr. Wakin received the DARPA Young Faculty Award for his research in compressive multi-signal processing for environments such as sensor and camera networks. In 2012, Dr. Wakin received the NSF CAREER Award for research into dimensionality reduction techniques for structured data sets. In 2014, Dr. Wakin received the Excellence in Research Award for his research as a junior faculty member at Mines. Dr. Wakin is a recipient of the Best Paper Award and the Signal Processing Magazine Best Paper Award from the IEEE Signal Processing Society. He has served as an Associate Editor for IEEE Signal Processing Letters and IEEE Transactions on Signal Processing, and he is currently a Senior Area Editor for IEEE Transactions on Signal Processing.
\end{IEEEbiography}

\appendices
%
%
%
%

\section{Proof of Theorem~\ref{thm:main_r1}}
\label{sec:proof_main_r1}

To prove Theorem~\ref{thm:main_r1}, we need the following lemma, which is proved in Appendix~\ref{proof_lem:xMxeig_r1}.
\begin{Lemma}\label{lem:xMxeig_r1}
Denote $\{\lambda_n\}_{n=1}^N$ as the eigenvalues of $\M$. Without loss of generality,  assume that $\lambda_1 > \lambda_2 \geq \lambda_3 \geq \cdots \geq \lambda_N$. Denote $\v_1$ as the eigenvector associated with $\lambda_1$. If $\|\text{{\em grad}}~g(\x)\|_2 \leq \epsilon$, then there exists some $n$ such that
\begin{align}
|\x^\top \M \x - \lambda_n| \leq \epsilon.
\label{eqn:xMxeig1}	
\end{align}	
Further, we have\footnote{Note that we assume $\langle\x,\v_1\rangle \geq 0$ here. In the case when $\langle\x,\v_1\rangle < 0$, one can bound $\|\x+\v_1\|_2^2$ instead.}
\begin{align}
\|\x-\v_1\|_2^2 & \leq \frac{2\epsilon}{\lambda_1-\lambda_2},  \quad \text{if}~n=1, \label{eqn:dist_xv1}\\
\|\x-\v_1\|_2^2 & \geq 2\left( 1-\frac{\epsilon}{\lambda_1-\lambda_n-\epsilon} \right), \quad  \text{if}~n\neq1.	\label{eqn:dist_xv2}
\end{align}
\end{Lemma}

Next, with the two bounds given in~\eqref{eqn:dist_xv1} and~\eqref{eqn:dist_xv2}, we show that there exist some positive constants $\epsilon$ and $\eta$ such that $|\lambda_{\min}(\text{hess}~g(\x))| \geq \eta$ when $\|\text{grad}~g(\x)\|_2 \leq \epsilon$. Let $\u\in\RRR^N$ be any vector that belongs to the tangent space of $\SSS^{N-1}$, i.e., $\u\in \TT_{\x}\SSS^{N-1} = \{\bxi \in \RRR^N: \x^\top \bxi = 0\}$. Without loss of generality, also let $\|\u\|_2=1$. Note that the quadratic term of the Riemannian Hessian is given as
\begin{align*}
\u^\top \text{hess}~g(\x) \u = \x^\top \M\x -\u^\top \M \u.	
\end{align*}

For $n=1$, we define $\ut = \V^\top \u$. Note that the square of the first entry in $\ut$ can be bounded with
\begin{align*}
\widetilde{u}_1^2 &= \langle \ut,\e_1\rangle^2 = \langle \u,\v_1\rangle^2
= \langle \u,\v_1-\x\rangle^2	\\
&\leq \|\u\|_2^2 \|\v_1-\x\|_2^2
\leq \frac{2\epsilon}{\lambda_1-\lambda_2}.
\end{align*}
Then, we can bound $\u^\top\M \u$ with
\begin{align*}
\u^\top\M \u &=\! \sum_{k=1}^N \lambda_k \widetilde{u}_k^2	
\!=\! \lambda_1 \widetilde{u}_1^2 + \sum_{k=2}^N \lambda_k \widetilde{u}_k^2
\leq \lambda_1 \widetilde{u}_1^2 + \lambda_2(1-\widetilde{u}_1^2) \\&= \lambda_2 + (\lambda_1 - \lambda_2) \widetilde{u}_1^2
\leq \lambda_2+2\epsilon.
\end{align*}
It follows that
\begin{align*}
\u^\top \text{hess}~g(\x) \u &\geq \lambda_1 - \epsilon - \lambda_2 - 2\epsilon
= \lambda_1 - \lambda_2 - 3\epsilon, 	
\end{align*}
which implies that $\lambda_{\min}(\text{hess}~g(\x)) \geq \eta$ as long as $3\epsilon+\eta \leq \lambda_1 - \lambda_2$.

For $n\neq 1$, with inequality~\eqref{eqn:dist_xv2}, we can get
\begin{align}
\langle \x,\v_1 \rangle &= 1-\frac 1 2 \|\x-\v_1\|_2^2
\leq \frac{\epsilon}{\lambda_1-\lambda_n-\epsilon}. 	
\label{eqn:innerxv}
\end{align}
By setting the direction $\u$ as $\u = \frac{(\I_N-\x\x^\top)\v_1}{\|(\I_N-\x\x^\top)\v_1\|_2}$ and noticing that $\|\u\|_2 = 1$ and $\u^\top\x=0$, we have
\begin{align*}
\u^\top \M \u &= \frac{\lambda_1 + (\x^\top \M \x - 2\lambda_1) \langle \x,\v_1 \rangle^2}{1-\langle \x,\v_1 \rangle^2}	\\
& \geq \frac{\lambda_1 - (2\lambda_1 - \lambda_n + \epsilon) \langle \x,\v_1 \rangle^2}{1-\langle \x,\v_1 \rangle^2}\\
& = (2\lambda_1 - \lambda_n + \epsilon) - \frac{\lambda_1 - \lambda_n + \epsilon}{1-\langle \x,\v_1 \rangle^2}\\
& \geq (2\lambda_1 - \lambda_n + \epsilon) - \frac{\lambda_1 - \lambda_n + \epsilon}{1-\frac{\epsilon^2}{(\lambda_1-\lambda_n-\epsilon)^2}},
\end{align*}
where the last inequality follows from~\eqref{eqn:innerxv}.
Then, the quadratic term of Riemannian Hessian can be bounded with
\begin{align}
&\u^\top \text{hess}~g(\x) \u \nonumber\\
\leq &\lambda_n + \epsilon -(2\lambda_1 - \lambda_n + \epsilon) + \frac{\lambda_1 - \lambda_n + \epsilon}{1-\frac{\epsilon^2}{(\lambda_1-\lambda_n-\epsilon)^2}} \nonumber \\
 = &-2(\lambda_1 - \lambda_n) + \frac{\lambda_1 - \lambda_n + \epsilon}{1-\frac{\epsilon^2}{(\lambda_1-\lambda_n-\epsilon)^2}}. \label{eqn:uhessu1}
\end{align}

Define a function $h(\alpha)$ as
\begin{align*}
h(\alpha) &\triangleq -2\alpha + \frac{\alpha + \epsilon}{1-\frac{\epsilon^2}{(\alpha-\epsilon)^2}}
=\frac{-\alpha^3+3\alpha^2 \epsilon - \alpha \epsilon^2 + \epsilon^3}{\alpha^2-2\alpha\epsilon}
\end{align*}
with $\epsilon >0$ being a fixed parameter and $\alpha > 2\epsilon$. With some fundamental calculations, we obtain the derivative of $h(\alpha)$, namely,
\begin{align*}
h'(\alpha) = \frac{-\alpha^4+4\alpha^3\epsilon - 5\alpha^2\epsilon^2-2\alpha \epsilon^3 + 2\epsilon^4}{(\alpha^2-2\alpha\epsilon)^2}.	
\end{align*}


Next, we argue that the above function $h(\alpha)$ decreases as we increase $\alpha$ by showing $h'(\alpha) < 0$ when $\alpha > 2\epsilon$. It is equivalent to show that
\begin{align*}
\hwt(\alpha) \triangleq \alpha^4-4\alpha^3\epsilon + 5\alpha^2\epsilon^2+2\alpha \epsilon^3 - 2\epsilon^4 > 0
\end{align*}
when $\alpha > 2\epsilon$. Note that
\begin{align}
\hwt'(\alpha)&= 4\alpha^3-12\alpha^2\epsilon+10\alpha\epsilon^2+2\epsilon^3, \nonumber\\
\hwt''(\alpha)&= 12\alpha^2-24\alpha\epsilon + 10\epsilon^2, \nonumber\\
\hwt'''(\alpha)&=24\alpha -24\epsilon. \label{eqn:dev3_ht}
\end{align}
Plugging in $\alpha = 2\epsilon$, we obtain
\begin{align}
\hwt(2\epsilon)&= 16\epsilon^4-32\epsilon^4+20\epsilon^4 + 4\epsilon^4-2\epsilon^4 = 6\epsilon^4>0, \nonumber\\	
\hwt'(2\epsilon)&= 32\epsilon^3-48\epsilon^3+20\epsilon^3+2\epsilon^3=6\epsilon^4>0, \nonumber\\
\hwt''(2\epsilon)&=48\epsilon^2 -48\epsilon^2+10\epsilon^2=10\epsilon^2>0. \label{eqn:dev2_ht_2ep}
\end{align}
It can be seen from~\eqref{eqn:dev3_ht} that $\hwt'''(\alpha) > 0$ for any $\alpha > 2\epsilon$. Together with~\eqref{eqn:dev2_ht_2ep}, we get
$\hwt''(\alpha) > \hwt''(2\epsilon) > 0, \quad \forall~\alpha > 2\epsilon,$	
 which further implies that
$\hwt'(\alpha) > \hwt'(2\epsilon) > 0, \quad \forall~\alpha > 2\epsilon,$	
and finally,
$\hwt(\alpha) > \hwt(2\epsilon) > 0, \quad \forall~\alpha > 2\epsilon.$
Therefore, we have $h'(\alpha)<0$ and thus $h(\alpha)$ is a monotonically decreasing function when $\alpha > 2\epsilon$. It follows from~\eqref{eqn:uhessu1} that
$\u^\top \text{hess}~g(\x) \u \leq h(\lambda_1-\lambda_n) \leq h(\lambda_1-\lambda_2)$
if $\lambda_1-\lambda_2 > 2\epsilon$. Then, we get
\begin{align*}
\u^\top \text{hess}~g(\x) \u
\leq & h(\lambda_1-\lambda_2)
\!=\!-2(\lambda_1 \!-\! \lambda_2) \!+\! \frac{\lambda_1 \!- \lambda_2 + \epsilon}{1-\frac{\epsilon^2}{(\lambda_1-\lambda_2-\epsilon)^2}}\\
= &  -0.72(\lambda_1-\lambda_2)
\leq -0.3(\lambda_1-\lambda_2)
\end{align*}
by setting $\epsilon = 0.2(\lambda_1-\lambda_2)$ since the value of second line increases as we increase $\epsilon$. Thus, we can choose $\eta = 0.3(\lambda_1-\lambda_2)$, which also satisfies $3\epsilon+\eta \leq \lambda_1 - \lambda_2$ required when $n=1$.
Therefore, we have verified that there do exist positive numbers $\epsilon$ and $\eta$ such that
$
|\lambda_{\min}(\text{hess}~g(\x))| \geq \eta	
$
when $\|\text{grad}~g(\x)\|_2 \leq \epsilon$.

\subsection{Proof of Lemma~\ref{lem:xMxeig_r1}}
\label{proof_lem:xMxeig_r1}

Denote $\M = \V \bLambda \V^\top$ as an eigendecomposition of $\M$. $\bLambda = \text{diag}([\lambda_1~\lambda_2~\cdots~\lambda_N])$ is a diagonal matrix that contains the eigenvalues of $\M$. $\V = [\v_1, \v_2, \cdots, \v_N]$ is an orthogonal matrix that contains the eigenvectors of $\M$.
Note that
\begin{align*}
&\|\text{grad}~g(\x)\|_2^2 \\
 = &\| (\x^\top \M \x)\x - \M\x \|_2^2  = \| (\x^\top \M \x)\x - \V \bLambda \V^\top\x \|_2^2\\
 \stack{\ding{172}}{=}& \| (\x^\top \M \x)\y - \bLambda \y \|_2^2
= \sumn (\x^\top \M\x-\lambda_n)^2 y_n^2\\
 \geq & \min_{n} (\x^\top \M\x-\lambda_n)^2 \sumn y_n^2
\stack{\ding{173}}{=} \min_{n} (\x^\top \M\x-\lambda_n)^2
\end{align*}
where \ding{172} follows by plugging in $\x = \V \y$ and from the fact that $\V$ is an orthogonal matrix, and \ding{173} follows from $\|\y\|_2^2 = \|\x\|_2^2 = 1$. Then, we have
$\min_{n}~ (\x^\top \M\x-\lambda_n)^2 \leq  \epsilon^2,$	
which implies that there exists some $n$ such that~\eqref{eqn:xMxeig1} holds.

Note that
\begin{align}
\|\x-\v_1\|_2^2 = \|\y-\e_1\|_2^2 = 2(1-y_1). \label{eqn:dist_xv0}	
\end{align}
Thus, bounding $\|\x-\v_1\|_2^2$ is equivalent to bounding $y_1$. We first consider the case when $n=1$, namely, $|\x^\top \M \x - \lambda_1| \leq \epsilon$. Note that
\begin{align*}
&\x^\top \M \x - \lambda_1 \stack{\ding{172}}{=} \sum_{k=1}^N \lambda_k y_k^2 - \lambda_1 = \lambda_1(y_1^2-1) + \sum_{k=2}^N \lambda_k y_k^2\\
\leq &\lambda_1(y_1^2-1) + \lambda_2 \sum_{k=2}^N  y_k^2
\stack{\ding{173}}{=} (\lambda_2-\lambda_1)(1-y_1)(1+y_1)
\end{align*}
where \ding{172} follows by plugging in $\M=\V\bLambda \V^\top$ and $\x = \V\y$, and \ding{173} follows from $\sum_{k=2}^N  y_k^2 = 1-y_1^2$. Then, we have
\begin{align*}
(\lambda_2-\lambda_1)(1-y_1)(1+y_1) \geq -\epsilon,	
\end{align*}
which implies that
\begin{align*}
1-y_1 \leq \frac{\epsilon}{(1+y_1)(\lambda_1-\lambda_2)}\leq \frac{\epsilon}{\lambda_1-\lambda_2}
\end{align*}
since $y_1 = \langle \x,\v_1\rangle \geq 0$. Plugging the above inequality into~\eqref{eqn:dist_xv0}, we can get~\eqref{eqn:dist_xv1}.

Next, we consider the case when $n\neq 1$. Observe that
\begin{align*}
\|\text{grad}~g(\x)\|_2^2  &= \sum_{k=1}^N (\x^\top \M\x-\lambda_k)^2 y_k^2
\geq 	(\x^\top \M\x-\lambda_1)^2 y_1^2,
\end{align*}
which further gives
$
(\x^\top \M\x-\lambda_1)^2 y_1^2 \leq  \epsilon^2.
$
It follows that $(\lambda_1-\x^\top \M\x)y_1 \leq \epsilon$. Together with $\x^\top \M \x \leq \lambda_n + \epsilon$, we have
\begin{align*}
y_1 \leq \frac{\epsilon}{\lambda_1-\x^\top \M\x}	 \leq \frac{\epsilon}{\lambda_1-\lambda_n - \epsilon}.
\end{align*}
Plugging the above inequality into~\eqref{eqn:dist_xv0}, we can then get~\eqref{eqn:dist_xv2} and finish the proof of Lemma~\ref{lem:xMxeig_r1}.

\section{Proof of Theorem~\ref{thm:main0}}
\label{sec:proof_main0}

For simplicity, we consider a diagonal matrix $\M$ in this proof, i.e., $\M = \diag([\lambda_1~\lambda_2~ \cdots~ \lambda_N])$ with $\lambda_1 > \lambda_2 >\cdots > \lambda_r > \lambda_{r+1} \geq \cdots \geq \lambda_N \geq 0$. For the case when $\M$ is not a diagonal matrix, one can diagonalize it with some unitary matrix. For the case when $\M$ has negative eigenvalues, one can add a constant to $\M$\footnote{Note that the eigendecomposition problem~\eqref{eqn:max_eigen} is equivalent to maximizing $\tr(\X^\top (\M+m_c\I) \X \N)$ on the Stiefel manifold. One can choose a positive constant $m_c$ such that $m_c + \lambda_N \geq 0$.} to ensure all the eigenvalues are non-negative.


Define $\Omega \triangleq \{i_1,\cdots i_r\}$ as a subset of $[N]$. Denote $\X_\Omega = [\e_{i_1},\cdots,\e_{i_r}] \in \RRR^{N\times r}$ with $\e_{i_j}$ being the $i_j$-th column of an identity matrix $\I_N$. Then, $\X_\Omega$ is a critical point of $g(\X)$. Note that
$
\X_\Omega^\top \M \X_\Omega = \text{diag}([
\lambda_{i_1} ~\cdots~ \lambda_{i_r}
])
\triangleq \bLambda_\Omega,
$
where the eigenvalues contained in $\bLambda_\Omega$ are not assumed to be in descending or ascending order.

For any $\U$ that belongs to the tangent space of $\St(N,r)$ at $\X_\Omega$, we have
\begin{align}
\X_\Omega^\top \U + \U^\top \X_\Omega =  \U(\Omega,:) + 	\U(\Omega,:)^\top = \zero,
\label{eq:UTU0}
\end{align}
where $\U(\Omega,:) \in \RRR^{r\times r}$ is a matrix that contains $r$ rows of $\U$ indexed by $\Omega$. Similarly, we use $\U(\Omega^c,:)\in \RRR^{(N-r)\times r}$ and $\bLambda_{\Omega^c} \in \RRR^{(N-r)\times (N-r)}$ to denote the  matrices that contain the remaining $N-r$ rows of $\U$ and $N-r$ eigenvalues of $\M$. Then, we have
\begin{equation}
\begin{aligned}
\U^\top \U &\!=\! \U(\Omega,:)^\top \U(\Omega,:) +  \U(\Omega^c,:)^\top \U(\Omega^c,:),\\
\U^\top \M \U & \!=\! \U(\Omega,:)\!^\top \!\bLambda_\Omega \! \U(\Omega,:) \!+\!  \U(\Omega^c,:)\!^\top \!\bLambda_{\Omega^c} \!\U(\Omega^c,:).	
\label{eq:UMU}
\end{aligned}
\end{equation}
Therefore, for any $\U\in \TT_{\X_\Omega}\St(N,r)$, we can rewrite the Riemannian Hessian as
\begin{align*}
&\text{hess}~g(\X_\Omega)[\U,\U] \\
=& \langle \X_\Omega^\top \M \X_\Omega, \U^\top \U \N\rangle - \langle \M,\U\N\U^\top \rangle \\	
=&\langle \bLambda_\Omega \N, \U^\top \U \rangle - \langle \U^\top\M \U,\N \rangle \\
=& \underbrace{\langle \bLambda_\Omega \N, \U(\Omega,:)^\top \U(\Omega,:) \rangle - \langle  \U(\Omega,:)^\top \bLambda_\Omega \U(\Omega,:),\N \rangle}_{\bPhi_1} \\
&+\! \underbrace{\langle \bLambda_\Omega \N, \U(\Omega^c,:)^\top \U(\Omega^c,:) \rangle \!-\! \langle  \U(\Omega^c,:)^\top \bLambda_{\Omega^c} \U(\Omega^c,:),\N \rangle}_{\bPhi_2},
\end{align*}
where the last equality follows by plugging~\eqref{eq:UMU}.

Let $\{\x_s\}_{s=1}^r$ denote the columns of a matrix $\X\in\RRR^{N\times r}$ and $\x_s(j)$ denote the $j$-th entry of $\x_s \in \RRR^{N}$. Note that
\begin{align*}
\bPhi_1 =& \langle \bLambda_\Omega \N, \U(\Omega,:)^\top \U(\Omega,:) \rangle - \langle  \U(\Omega,:)^\top \bLambda_\Omega \U(\Omega,:),\N \rangle \\
=& \sum_{s=1}^r \mu_s \lambda_{i_s} \sum_{j=1}^r \u_s^2(i_j)	 - \sum_{s=1}^r \mu_s \sum_{j=1}^r \lambda_{i_j}\u_s^2(i_j)\\
=& \sum_{s,j=1}^r \mu_s(\lambda_{i_s} - \lambda_{i_j} )\u_s^2(i_j)
\stack{\ding{172}}{=} \sum_{j,s=1}^r \mu_j(\lambda_{i_j} - \lambda_{i_s} )\u_j^2(i_s)\\
=& \frac 1 2 \sum_{s,j=1}^r \mu_s(\lambda_{i_s} \!-\! \lambda_{i_j} )\u_s^2(i_j) \!+\! \frac 1 2 \sum_{j,s=1}^r \mu_j(\lambda_{i_j} \!-\! \lambda_{i_s} )\u_j^2(i_s)\\
\stack{\ding{173}}{=}& \frac 1 2 \sum_{s,j=1}^r (\mu_s - \mu_j) (\lambda_{i_s} - \lambda_{i_j} )\u_s^2(i_j)\\
\stack{\ding{174}}{=}& \frac 1 2 \sum_{s,j=1 \atop s\neq j}^r (\mu_s - \mu_j) (\lambda_{i_s} - \lambda_{i_j} )\u_s^2(i_j)\\
=& \sum_{s,j=1 \atop s<j}^r (\mu_s - \mu_j) (\lambda_{i_s} - \lambda_{i_j} )\u_s^2(i_j).
\end{align*}
Here, \ding{172} follows by exchanging the role of indices $j$ and $s$. \ding{173} and \ding{174} follow from the fact that~\eqref{eq:UTU0} implies
\begin{align}
\u_s(i_j) = \begin{cases}
-\u_j(i_s),~~& s \neq j,\\
0, & s=j, 	
\end{cases}	
\quad \forall~j,s\in[r].
\label{eq:tangent_cond}
\end{align}
Similarly, we have
\begin{align*}
\bPhi_2 =& \langle \bLambda_\Omega \N, \U(\Omega^c,:)\!^\top\! \U(\Omega^c,:) \rangle \!-\! \langle  \U(\Omega^c,:)\!^\top\! \bLambda_{\Omega^c} \U(\Omega^c,:),\N \rangle\\
=& \sum_{s=1}^r \mu_s \lambda_{i_s} 	\sum_{j\in\Omega^c} \u_s^2(j) - \sum_{s=1}^r \mu_s \sum_{j\in\Omega^c} \lambda_j \u_s^2(j)\\
=& \sum_{s=1}^r \sum_{j\in\Omega^c} \mu_s (\lambda_{i_s} - \lambda_j)\u_s^2(j).
\end{align*}

Next, we bound $\text{hess}~g(\X_\Omega)[\U,\U]$ for any $\U\in \TT_{\X_\Omega}\St(N,r)$ by considering the following three cases:
\begin{itemize}
\item {\bf Case 1:} $\Omega \triangleq \{i_1,\cdots i_r\} = [r]$.	
\item {\bf Case 2:} $\Omega \!\triangleq \!\{i_1,\cdots i_r\} \!=\! \PP_{mu}([r]) \!\neq \![r]$.
\item {\bf Case 3:} $\Omega \triangleq \{i_1,\cdots i_r\} \neq \PP_{mu}([r])$.
\end{itemize}

\subsection{Case 1}
Note that the two terms $\bPhi_1$ and $\bPhi_2$ can be bounded with
\begin{align*}
\bPhi_1 =& \sum_{s,j=1 \atop s<j}^r (\mu_s - \mu_j) (\lambda_{i_s} - \lambda_{i_j} )\u_s^2(i_j)
\stack{\ding{172}}{\geq}  \sum_{s,j=1 \atop s<j}^r (\lambda_s - \lambda_j )\u_s^2(j) \\
\geq &d_{\min} \sum_{s,j=1 \atop s<j}^r\u_s^2(j) = \frac 1 2 d_{\min} \|\U(\Omega,:)\|_F^2
\end{align*}
and
\begin{align*}
\bPhi_2 \!=\!& \sum_{s=1}^r \!\sum_{j\in\Omega^c} \mu_s (\lambda_{i_s} \!-\! \lambda_j)\u_s^2(j)	
\!\stack{\ding{173}}{=}\! \sum_{s=1}^r \!\sum_{j=r+1}^N \mu_s (\lambda_s \!-\! \lambda_j)\u_s^2(j) \\
\stack{\ding{174}}{\geq}& (\lambda_r-\lambda_{r+1})\|\U(\Omega^c,:)\|_F^2
\geq  d_{\min} \|\U(\Omega^c,:)\|_F^2
\end{align*}
where \ding{172} follows from $\mu_s - \mu_j \geq 1$, $i_s = s$ and $i_j=j$. \ding{173} follows from $\Omega^c = \{r+1,\cdots,N\}$. \ding{174} follows from $\mu_s \geq 1$. Then, we get
\begin{align*}
\text{hess}~g(\X_\Omega)[\U,\U] &= \bPhi_1 + \bPhi_2 \\ &\geq \frac 1 2 d_{\min} \|\U(\Omega,:)\|_F^2 +	 d_{\min} \|\U(\Omega^c,:)\|_F^2 \\
&\geq \frac 1 2 d_{\min} \|\U\|_F^2, \quad \forall~\U\in \TT_{\X_\Omega}\St(N,r),
\end{align*}
which implies that
\begin{align}
\lambda_{\min}(\text{hess}~g(\X_\Omega)) \geq \frac 1 2 d_{\min}.
\label{eq:localmin}	
\end{align}
Therefore, $\X_\Omega = [\e_1,\cdots,\e_r] \in \RRR^{N\times r}$ is a global minimum of $g(\X)$.

\subsection{Case 2}

As $\U$ belongs to the tangent space of $\St(N,r)$ at $\X_\Omega$, we need to construct a direction $\U$ that satisfies the condition~\eqref{eq:tangent_cond}. Now, we construct a direction $\U$ with $\U(\Omega^c,:) = \zero$. Then, we have $\bPhi_2 = 0$. Furthermore, for any $s,j\in[r]$ with $s<j$, we set
\begin{align*}
\u_s(i_j) = -\u_j(i_s) \begin{cases}
 = 0, \quad &\text{if}~i_s < i_j,\\
\neq 0, &\text{if}~i_s > i_j	.
\end{cases}	
\end{align*}
Then, we have
\begin{align*}
\bPhi_1 =& \sum_{s,j=1 \atop s<j}^r (\mu_s - \mu_j) (\lambda_{i_s} - \lambda_{i_j} )\u_s^2(i_j)\\
\stack{\ding{172}}{=}&	-\sum_{s,j=1 \atop s<j,i_s > i_j}^r (\mu_s - \mu_j) (\lambda_{i_j} - \lambda_{i_s} )\u_s^2(i_j)\\
\stack{\ding{173}}{\leq}& -d_{\min} \sum_{s,j=1 \atop s<j,i_s > i_j}^r \u_s^2(i_j)
\stack{\ding{174}}{=} - \frac 1 2 d_{\min} \|\U\|_F^2,
\end{align*}
where \ding{172} follows from $\u_s(i_j) = 0$ if $i_s < i_j$. \ding{173} follows from $\mu_s-\mu_j \geq 1$ and
$\lambda_{i_j} - \lambda_{i_s} \geq \min_{1\leq i_j < i_s \leq r} \lambda_{i_j} - \lambda_{i_s} \geq d_{\min}.$
\ding{174} follows from
\begin{align*}
\sum_{s,j=1 \atop s<j,i_s > i_j}^r \!\!\!\!\!\!\!\u_s^2(i_j) \!= \!\!\sum_{s,j=1 \atop s<j}^r \!\!\!\u_s^2(i_j)	 \!=\!\!\! \sum_{s,j=1 \atop j<s}^r \!\!\!\u_j^2(i_s) \!= \!\!\frac 1 2	\!\!\sum_{s,j=1 \atop s\neq j}^r\!\! \u_s^2(i_j)	 \!=\! \frac 1 2 \!\|\U\|_F^2.
\end{align*}

It follows that there exists some $\U\!\in \!\TT_{\X_\Omega}\St(N,r)$ such that
\begin{align*}
\text{hess}~g(\X_\Omega)[\U,\U] &= \bPhi_1 + \bPhi_2 \leq - \frac 1 2 d_{\min} \|\U\|_F^2,
\end{align*}
which further implies
$\lambda_{\min}(\text{hess}~g(\X_\Omega)) \leq -\frac 1 2 d_{\min}.$	
Therefore,
$\X_\Omega = [\e_{i_1},\cdots,\e_{i_r}] \in \RRR^{N\times r}$ with $\Omega = \PP_{mu}([r]) \neq [r]$ are strict saddles of $g(\X)$.

\subsection{Case 3}

Note that there exist some $i^\star,j^\star \in[r]$ such that $i^\star \notin \Omega$ and $i_{j^\star} \in \Omega$ with $i_{j^\star} \geq r+1$. Thus, we have
$\lambda_{i^\star} - \lambda_{i_{j^\star}} \geq \lambda_r - \lambda_{r+1} \geq d_{\min}.$	
Now, we set $\U$ as a matrix with only one non-zero entry at $(i^\star,j^\star)$, namely, $\U(\Omega,:) = \zero$ and $\u_{j^\star}^2(i^\star) = \|\U\|_F^2$. It can be seen that such a direction $\U$ belongs to the tangent space of $\St(N,r)$ at $\X_\Omega$. Then, we have $\bPhi_1 = 0$ since $\U(\Omega,:) = \zero$. We can bound $\bPhi_2$ with
\begin{align*}
\bPhi_2 =&  \sum_{s=1}^r \sum_{j\in\Omega^c} \mu_s (\lambda_{i_s} - \lambda_j)\u_s^2(j)	\\
=& -\mu_{j^\star} (\lambda_{i^\star} - \lambda_{i_{j^\star}})\|\U\|_F^2	
\leq  -d_{\min} \|\U\|_F^2,
\end{align*}
which further implies that
\begin{align*}
\text{hess}~g(\X_\Omega)[\U,\U] &= \bPhi_1 + \bPhi_2 \leq -  d_{\min} \|\U\|_F^2.
\end{align*}
Finally, we have
$\lambda_{\min}(\text{hess}~g(\X_\Omega)) \leq - d_{\min}.	$
Therefore,
$\X_\Omega = [\e_{i_1},\cdots,\e_{i_r}] \in \RRR^{N\times r}$ with $\Omega \neq \PP_{mu}([r])$ are also strict saddles of $g(\X)$.

\section{Proof of Theorem~\ref{thm:main}}
\label{sec:proof_main}

As in Appendix~\ref{sec:proof_main0}, we consider a diagonal matrix $\M$ to simplify the proof, i.e., $\M = \diag([\lambda_1~\lambda_2~ \cdots~ \lambda_N])$ with $\lambda_1 > \lambda_2 >\cdots > \lambda_r > \lambda_{r+1} \geq \cdots \geq \lambda_N \geq 0$.

\begin{Lemma}\label{lem:xMxeig}
Denote $\{\lambda_n\}_{n=1}^N$ as the eigenvalues of $\M$. Without loss of generality,  assume that $\lambda_1 > \lambda_2 >\cdots \lambda_r > \lambda_{r+1} \geq \cdots \geq \lambda_N \geq 0$. Define an index set $\Omega \triangleq \{i_1,\cdots i_r\}$ as a subset of $[N]$. Denote $\bLambda_\Omega = \diag([\lambda_{i_1},\cdots,\lambda_{i_r}])\in\RRR^{r \times r}$ and $\X_\Omega = [\e_{i_1},\cdots,\e_{i_r}] \in \RRR^{N\times r}$ as a diagonal matrix that contains $r$ eigenvalues of $\M$ and a matrix that contains the corresponding eigenvectors.
If $\|\text{{\em grad}}~g(\X)\|_F\leq \epsilon$, then there exists an index set $\Omega$ such that
\begin{align}
\|\X^\top \M \X - \bLambda_\Omega\|_F \leq 4\epsilon.
\label{eq:xMxeig1}	
\end{align}	
Moreover, we have
\begin{align}
\|\X-\X_{\Omega}\|_F &\leq 12 d_{\min}^{-1} \epsilon, \quad \text{if}~~ \Omega = [r].
\label{eq:dist_xv1}
\end{align}
\end{Lemma}
The above Lemma is proved in Appendix~\ref{proof_lem:xMxeig}.
Next, we bound $\text{hess}~g(\X)[\U,\U]$
in
the following three cases:

\begin{itemize}
\item {\bf Case 1:} $\Omega \triangleq \{i_1,\cdots i_r\} = [r]$.
\item {\bf Case 2:} $\Omega \!\triangleq \!\{i_1,\cdots i_r\}\! =\! \PP_{mu}([r]) \!\neq \![r]$.
\item {\bf Case 3:} $\Omega \triangleq \{i_1,\cdots i_r\} \neq \PP_{mu}([r])$.
\end{itemize}

\subsection{Case 1}

For any $\U\!\in\! \TT_{\X}\St(N,r)$, we bound $\text{hess}~g(\X)[\U,\!\U]$ as
\begin{align}
\text{hess}~g(\X)&[\U,\U]
\geq \underbrace{\text{hess}~g(\X_\Omega)[\U,\U]}_{\textbf{Term I}} \nonumber \\ -& \underbrace{| \text{hess}~g(\X)[\U,\U] - \text{hess}~g(\X_\Omega)[\U,\U] |}_{\textbf{Term II}}.
\label{eq:hess_case1_1}
\end{align}
Recall that
\begin{align*}
\text{hess}~g(\X)[\U,\U] &=  \langle \X^\top \M \X, \U^\top \U \N\rangle - \langle \M,\U\N\U^\top \rangle\\
\text{hess}~g(\X_\Omega)[\U,\U] &=  \langle \bLambda_\Omega, \U^\top \U \N\rangle - \langle \M,\U\N\U^\top \rangle.
\end{align*}
Then, we have
\begin{equation}
\begin{aligned}
\textbf{Term II} =& | \langle ( \X^\top \M \X - \bLambda_\Omega )\N,\U^\top \U \rangle |	\\
\leq & \|\N\| \|\X^\top \M \X - \bLambda_\Omega\|_F \|\U\|_F^2
\leq  4r\epsilon \|\U\|_F^2,
\label{eq:hess_case1_2}
\end{aligned}
\end{equation}
where the last inequality follows from $\|\N\| = r$ and~\eqref{eq:xMxeig1}.

Denote $\PP_{\TT_{\X}}$ as orthogonal projection onto $\TT_{\X}\St(N,r)$. Then, we have
\begin{equation}
\begin{aligned}
\textbf{Term I} =& \text{hess}~g(\X_\Omega)[\PP_{\TT_{\X_\Omega}}\U,\PP_{\TT_{\X_\Omega}}\U]	\\
\geq &  \lambda_{\min} (\text{hess}~g(\X_\Omega)) \| \PP_{\TT_{\X_\Omega}}\U \|_F^2\\
\geq & \frac 1 2 d_{\min} ( 1-2\|\X-\X_\Omega\|_F )^2\|\U\|_F^2,
\label{eq:hess_case1_3}
\end{aligned}
\end{equation}
where the last inequality follows from~\eqref{eq:localmin} and
\begin{align*}
&\| \PP_{\TT_{\X_\Omega}}\!\! \U \|_F \!=\! \| \U \!+\! \PP_{\TT_{\X_\Omega}}\!\!\U \!-\! \U \|_F
\!= \!	\| \U \!+\! \PP_{\TT_{\X_\Omega}}\!\!\U \!-\! \PP_{\TT_{\X}}\!\U \|_F\\
\geq & \|\U\|_F - \|\PP_{\TT_{\X_\Omega}}\U - \PP_{\TT_{\X}}\U \|_F\\
\stack{\ding{172}}{=}& \|\U\|_F - \frac 1 2\|\X(\X^\top \U + \U^\top \X) - \X_\Omega(\X_\Omega^\top \U + \U^\top \X_\Omega) \|_F\\
=& \|\U\|_F \!-\! \frac 1 2\|\X\X^\top \U - \X_\Omega\X_\Omega^\top \U + \X\U^\top \X  - \X_\Omega\U^\top \X_\Omega \|_F\\
\geq & \|\U\|_F \!\!-\! \frac 1 2\|\X\X^\top \!\U \!- \!\X_\Omega\X_\Omega^\top\! \U\|_F \!\!- \!\frac 1 2 \| \X\U^\top \!\X  \!-\!\X_\Omega\!\U^\top\! \X_\Omega \|_F\\
\stack{\ding{173}}{\geq} & (1-2\|\X-\X_\Omega\|_F) \|\U\|_F.
\end{align*}
Here, \ding{172} follows from $\PP_{\TT_{\X}}\U = \U - \frac 1 2 \X(\X^\top \U + \U^\top \X)$~\cite[Section 3.6.1]{absil2009optimization}. \ding{173} follows from
\begin{align*}
&\|\X\X^\top \U - \X_\Omega\X_\Omega^\top \U\|_F 	\\
=& \|\X\X^\top \U - \X\X_\Omega^\top \U + \X\X_\Omega^\top \U - \X_\Omega\X_\Omega^\top \U\|_F\\
\leq & \|\X\X^\top \U - \X\X_\Omega^\top \U\|_F + \|\X\X_\Omega^\top \U - \X_\Omega\X_\Omega^\top \U\|_F\\
= & \|\X(\X-\X_\Omega)^\top \U\|_F + \|(\X - \X_\Omega)\X_\Omega^\top \U\|_F\\
\leq & 2\|\X-\X_\Omega\|_F\|\U\|_F,
\\
\text{and} &
\\
&\| \X\U^\top \X  -\X_\Omega\U^\top \X_\Omega \|_F \\
=&	 \| \X\U^\top \X - \X\U^\top \X_\Omega + \X\U^\top \X_\Omega  -\X_\Omega\U^\top \X_\Omega \|_F\\
\leq & \| \X\U^\top \X - \X\U^\top \X_\Omega \|_F + \|\X\U^\top \X_\Omega  -\X_\Omega\U^\top \X_\Omega \|_F\\
= & \| \X\U^\top (\X - \X_\Omega) \|_F + \|(\X  -\X_\Omega)\U^\top \X_\Omega \|_F\\
\leq & 2\|\X-\X_\Omega\|_F\|\U\|_F.
\end{align*}

It remains to bound $\|\X-\X_\Omega\|_F$ with $\Omega = [r]$. It follows from~\eqref{eq:dist_xv1}-\eqref{eq:hess_case1_3}  that
\begin{align*}
&\text{hess}~g(\X)[\U,\U] \!\geq \! 	\left[\frac 1 2 d_{\min} ( 1\!-\!2\|\X\!-\!\X_\Omega\|_F )^2 \!-\! 4r\epsilon \right]\!\|\U\|_F^2  \\
\geq & \left[\frac 1 2 d_{\min} ( 1-4\|\X-\X_\Omega\|_F ) - 4r\epsilon \right]\|\U\|_F^2  \\
\geq & \left[\!\frac 1 2 d_{\min} ( 1\!-\!48  d_{\min}^{-1}\epsilon ) \!-\! 4r\epsilon \right]\!\!\|\U\|_F^2
\!=\!  \left(\frac 1 2 d_{\min} \!-\! 24 \epsilon \!-\! 4r\epsilon\right)\!\!\|\U\|_F^2.
\end{align*}
Finally, we get
$\lambda_{\min}(\text{hess}~g(\X)) \geq  \frac 1 2 d_{\min} - 24  \epsilon - 4r\epsilon.$

\subsection{Case 2}

In this case, we will show that there exists a $\U \in \TT_{\X}\St(N,r)$ such that $\text{hess}~g(\X)[\U,\U] <0$. Recall that
\begin{align*}
 \TT_{\X}\St(N,r) =	\{\X \S+ \X_\perp \K: \S^\top = -\S, \K \in \RRR^{(N-r) \times r} \}.
\end{align*}
For any $s,j\in[r]$ with $s<j$ and $i_s > i_j$, we construct $\U$ by setting
$\K = \zero, \quad \text{and} \quad \S= \e_s\e_j^\top - \e_j \e_s^\top.$	

We can then bound the Riemannian Hessian as
\begin{align*}
&\text{hess}~g(\X)[\U,\U] \\
=&  \langle \X^\top \M \X \N, \U^\top \U \rangle - \langle \U^\top\M\U,\N \rangle\\
=& \langle \X^\top \M \X \N, \S^\top \S \rangle - \langle \S^\top \X^\top \M \X \S,\N \rangle \\
=& \langle \X^\top \M \X \N, \e_j\e_j^\top + \e_s \e_s^\top \rangle\\
&~~~~~~ - \langle (\e_s\e_j^\top - \e_j \e_s^\top)^\top \X^\top \M \X (\e_s\e_j^\top - \e_j \e_s^\top),\N \rangle\\
=& \mu_j \x_j^\top \M \x_j + \mu_s \x_s^\top \M \x_s - \mu_j\x_s^\top \M \x_s - \mu_s \x_j^\top \M \x_j \\
=& (\mu_s - \mu_j) (\x_s^\top \M \x_s - \x_j^\top \M \x_j)\\
=& (\mu_s \!-\! \mu_j)\! \left[(\x_s^\top \M \x_s \!-\! \lambda_{i_s}) \!-\! ( \x_j^\top \M \x_j \!-\! \lambda_{i_j}) \!+\! (\lambda_{i_s} \!-\! \lambda_{i_j})  \right]\\
\leq & (3\epsilon - \frac 1 2 d_{\min})\|\U\|_F^2,
\end{align*}
where the last inequality follows from $\mu_s - \mu_j \geq 1$, $\lambda_{i_s} - \lambda_{i_j} \leq -d_{\min}$, $\|\U\|_F^2 = \|\S\|_F^2 = 2$, and the following inequality
\begin{align}
\lambda_{i_j} - 3\epsilon \leq  \x_j^\top\M \x_j \leq \lambda_{i_j} + 3\epsilon, \quad \forall~j\in[r],
\label{eq:xMx3e}
\end{align}
which can be obtained by combining inequalities~\eqref{eq:xmxl1} and~\eqref{eq:xmxl2}. In particular, one has
\begin{align*}
( \x_j^\top\M \x_j - \lambda_{i_j} )^2 \leq 	\|\X\D_{\M} - \M\X\|_F ^2 \leq 9\epsilon^2, \quad \forall~j\in[r],
\end{align*}
which further indicates the inequality~\eqref{eq:xMx3e}.

Then, we have
$\lambda_{\min}(\text{hess}~g(\X)) \leq -\frac 1 2 d_{\min} + 3\epsilon.$

\subsection{Case 3}

Note that there exist some $i^\star,j^\star \in [r]$ such that $i^\star \notin \Omega$ and $i_{j^\star} \in \Omega$ with $i_{j^\star} \geq r+1$. Then, we have
$\lambda_{i^\star} \geq \lambda_r \quad \text{and} \quad \lambda_{i_{j^\star}} \leq \lambda_{r+1}.	$
Next, we construct $\U$ by setting
$\S = \zero, \quad \text{and} \quad \K= \X_{\perp}^\top \e_{i^\star} \e_{j^\star}^\top,	$
where $\e_{i^\star} \in \RRR^N$ and $\e_{j^\star}\in \RRR^r$ are the $i^\star$-th column of an identity matrix $\I_N$ and the $j^\star$-th column of an identity matrix $\I_r$. Note that $\|\U\|_F^2 = \|\X_\perp \K \|_F^2 \leq 1$.

Then, we can bound the Riemannian Hessian as
\begin{equation}
\begin{aligned}
&\text{hess}~g(\X)[\U,\U] =  \langle \X^\top \M \X \N, \U^\top \U \rangle - \langle \U^\top\M\U,\N \rangle\\
=& \e_{i^\star}^\top \X_\perp \X_\perp^\top \e_{i^\star} \cdot \mu_{j^\star} \x_{j^\star}^\top \M \x_{j^\star} \!-\! \e_{i^\star}^\top \X_\perp \X_\perp^\top \M \X_\perp \X_\perp^\top \e_{i^\star} \cdot \mu_{j^\star}\\
=& \mu_{j^\star} \|\X_\perp^\top \e_{i^\star}\|_2^2 \left( \x_{j^\star}^\top \M \x_{j^\star} - \frac{\e_{i^\star}^\top \X_\perp \X_\perp^\top \M \X_\perp \X_\perp^\top \e_{i^\star} }{ \|\X_\perp^\top \e_{i^\star}\|_2^2  } \right),
\label{eq:hesscase30}
\end{aligned}
\end{equation}
where the second equality follows by plugging $\U = \X_\perp \K = \X_\perp \X_\perp^\top \e_{i^\star} \e_{j^\star}^\top$.
Note that
\begin{equation}
\begin{aligned}
&\e_{i^\star}^\top \X_\perp \X_\perp^\top \M \X_\perp \X_\perp^\top \e_{i^\star}  \!=\! (\e_{i^\star} \!\!-\! \X\X^\top\! \e_{i^\star}\!)\!^\top \!\M (\e_{i^\star} \!-\! \X\X^\top \!\e_{i^\star}\!) 	\\
=& \lambda_{i^\star} - 2\lambda_{i^\star} \|\X^\top \e_{i^\star} \|_2^2 + \e_{i^\star}^\top \X \X^\top \M \X \X^\top \e_{i^\star}.
\label{eq:expxpMxpxpe}
\end{aligned}	
\end{equation}
Recall that we denote $\X^\top \M \X = \D_{\M} + \E_{\M}$
with $\D_{\M}$ and $\E_{\M}$ being the diagonal and off-diagonal parts of $\M$ in equation~\eqref{eq:xmxde}. Then, one can bound
\begin{equation}
\begin{aligned}
& \e_{i^\star}^\top \X \X^\top \M \X \X^\top \e_{i^\star} =  \e_{i^\star}^\top \X \D_{\M} \X^\top \e_{i^\star}+ \e_{i^\star}^\top \X\E_{\M} \X^\top \e_{i^\star}	\\
=& \e_{i^\star}^\top\M \X\X^\top \!\e_{i^\star} \!+\! \e_{i^\star}^\top (\X \D_{\M} \!-\! \M \X) \X^\top \!\e_{i^\star}\!+ \! \e_{i^\star}^\top \X\E_{\M} \X^\top \!\e_{i^\star},\\
\stack{\ding{172}}{\geq} & \lambda_{i^\star} \!\| \X\!^\top \!\e_{i^\star}\!\|_2^2 \!\!-\!\! \| \e_{i^\star}\! \|_2 \|\X \D_{\M} \!\!-\!\! \M \X\|\!_F \|\X\!^\top\! \e_{i^\star}\!\|_2 \!\!-\!\! \|\E_{\M}\|\!_F \| \X\!^\top\! \!\e_{i^\star}\!\|_2^2\\
\stack{\ding{173}}{\geq} & \lambda_{i^\star} \| \X^\top \e_{i^\star}\|_2^2 - 5\epsilon,
\label{eq:exxmxxe}
\end{aligned}	
\end{equation}
where \ding{172} follows from $\e_{i^\star}^\top\M \X\X^\top \e_{i^\star} = \lambda_{i^\star} \| \X^\top \e_{i^\star}\|_2^2$ and the Cauchy-Schwarz inequality. \ding{173} follows from inequalities~\eqref{eq:xmxl1} and~\eqref{eq:embound}.
Plugging~\eqref{eq:exxmxxe} into~\eqref{eq:expxpMxpxpe}, we obtain
\begin{align*}
\e_{i^\star}^\top \X_\perp \X_\perp^\top \M \X_\perp \X_\perp^\top \e_{i^\star} \geq & \lambda_{i^\star} - \lambda_{i^\star} \|\X^\top \e_{i^\star} \|_2^2 - 5\epsilon\\
=& \lambda_{i^\star} \|\X_\perp^\top \e_{i^\star} \|_2^2 - 5\epsilon,
\end{align*}	
which together with~\eqref{eq:hesscase30} gives
\begin{align*}
&\text{hess}~g(\X)[\U,\U]\\
\leq & \mu_{j^\star} \|\X_\perp^\top \e_{i^\star}\|_2^2 \left( \x_{j^\star}^\top \M \x_{j^\star} - \lambda_{i^\star} +  \frac{ 5\epsilon}{ \|\X_\perp^\top \e_{i^\star}\|_2^2  } \right)\\
\stack{\ding{172}}{\leq} & \mu_{j^\star} \|\X_\perp^\top \e_{i^\star}\|_2^2 \left( \lambda_{r+1} + 3\epsilon - \lambda_r +  \frac{ 5\epsilon}{ \|\X_\perp^\top \e_{i^\star}\|_2^2  } \right)\\
\stack{\ding{173}}{\leq} &  \mu_{j^\star} \|\X_\perp^\top \e_{i^\star}\|_2^2 \left( -d_{\min} +   \frac{ 8\epsilon}{ \|\X_\perp^\top \e_{i^\star}\|_2^2  } \right)\|\U\|_F^2\\
\stack{\ding{174}}{\leq} & \mu_{j^\star} \|\X_\perp^\top \e_{i^\star}\|_2^2 \left( -d_{\min} +   \frac{ 8\epsilon}{1 - 36 d_{\min}^{-2}	\epsilon^2} \right) \|\U\|_F^2
\end{align*}
Here, \ding{172} follows from $\lambda_{i^\star} \geq \lambda_r$ and $\x_{j^\star}^\top\M \x_{j^\star} \leq \lambda_{i_{j^\star}} + 3\epsilon \leq \lambda_{r+1} + 3\epsilon$.
\ding{173} follows from $\lambda_r-\lambda_{r+1} \geq d_{\min}$, $\|\X_\perp^\top \e_{i^\star}\|_2^2 \leq 1$ and $\|\U\|_F^2 \leq 1$. \ding{174} follows from
\begin{align}
\|\X_\perp^\top \e_{i^\star}\|_2^2 \geq 	1 - 36 d_{\min}^{-2}	\epsilon^2.
\label{eq:xpe}
\end{align}
Then, we have
\begin{align*}
\lambda_{\min}(\text{hess}~g(\X)) \leq \mu_{j^\star} \|\X_\perp^\top \e_{i^\star}\|_2^2 \left( -d_{\min} +   \frac{ 8\epsilon}{ 	1 - 36 d_{\min}^{-2}	\epsilon^2 } \right).
\end{align*}

Finally, it remains to show that~\eqref{eq:xpe} holds in Case 3.
It follows from~\eqref{eq:xlomx} that
$\|\X\bLambda_{\Omega} - \M\X\|_F \leq  6\epsilon,$
which allows us to bound the $\ell_2$-norm of the $i^\star$-th row of $\X\bLambda_{\Omega} - \M\X$ as
\begin{align}
\| \X(i^\star,:)\odot [ \lambda_{i_1} ~\cdots ~ \lambda_{i_r} ] - \X(i^\star,:)\odot [ \lambda_{i^\star} ~\cdots ~ \lambda_{i^\star} ]    \|_2 \leq 6\epsilon,
\label{eq:rownorm1}	
\end{align}
where $\odot$ denotes the elementary-wise multiplication. On the other hand, we have
\begin{equation}
\begin{aligned}
&\| \X(i^\star,:)\odot [ \lambda_{i_1} ~\cdots ~ \lambda_{i_r} ] - \X(i^\star,:)\odot [ \lambda_{i^\star} ~\cdots ~ \lambda_{i^\star} ]    \|_2 \\
= & \| \X(i^\star,:)\odot [ \lambda_{i_1} -\lambda_{i^\star} ~\cdots ~ \lambda_{i_r} -\lambda_{i^\star} ]    \|_2	\\
\geq & \min_{j\in[r]} |\lambda_{i_j} -\lambda_{i^\star}| \cdot \| \X(i^\star,:)\|_2
\geq  d_{\min}\| \X(i^\star,:)\|_2.
\label{eq:rownorm2}
\end{aligned}
\end{equation}
Combining~\eqref{eq:rownorm1} and \eqref{eq:rownorm2} leads to
\begin{align*}
\| \X(i^\star,:)\|_2 \leq 6d_{\min}^{-1}	\epsilon,
\end{align*}
which immediately implies that
\begin{align*}
\|\X_\perp^\top \e_{i^\star}\|_2^2 = 1 \!-\! \|\X^\top \e_{i^\star}\|_2^2 = 1 \!-\! \| \X(i^\star,:) \|_2^2 \geq 1 \!-\! 36 d_{\min}^{-2}	\epsilon^2.
\end{align*}

\subsection{Summary}

\begin{itemize}
\item {\bf Case 1:} $\Omega \triangleq \{i_1,\cdots i_r\} = [r]$.
\begin{align*}
\lambda_{\min}(\text{hess}~g(\X)) \geq  \frac 1 2 d_{\min} - 24 \epsilon - 4r\epsilon \geq 0.11 d_{\min}
\end{align*}
if $\epsilon \leq \frac{1}{72} r^{-1} d_{\min}$.

\item {\bf Case 2:} $\Omega \triangleq \{i_1,\cdots i_r\} = \PP_{mu}([r]) \neq [r]$.
\begin{align*}
\lambda_{\min}(\text{hess}~g(\X)) \leq -\frac 1 2 d_{\min} + 3\epsilon \leq -0.45 d_{\min}
\end{align*}
if $\epsilon \leq \frac{1}{72} r^{-1} d_{\min}$.

\item {\bf Case 3:} $\Omega \triangleq \{i_1,\cdots i_r\} \neq \PP_{mu}([r])$.
\begin{align*}
\lambda_{\min}(\text{hess}~g(\X)) &\leq \mu_{j^\star}\! \|\X_\perp^\top \e_{i^\star}\!\|_2^2 \!\left(\!\! -d_{\min} \!+ \!  \frac{ 8\epsilon}{ 	1 \!-\! 36 d_{\min}^{-2}	\epsilon^2 } \!\!\right)\\
& \leq  -0.88 d_{\min}
\end{align*}
if $\epsilon \leq \frac{1}{72} r^{-1} d_{\min}$.

\end{itemize}

Therefore, we can set $\epsilon = \frac{1}{72} r^{-1} d_{\min}$ and $\eta = 0.11 d_{\min}$.

\section{Proof of Lemma~\ref{lem:xMxeig}}
\label{proof_lem:xMxeig}

Note that
\begin{align*}
\|\text{grad}~g(\X)\|_F^2 \!=\!& \left\|(\X\X^\top \!-\! \I_N)\M\X \N	\!-\!  \frac 1 2 \X[\X^\top \M \X, \N]\right\|_F^2 \\
\!=& \|(\X\X^\top \!\!-\! \I_N)\M\X \N\|_F^2	 \!+\! \frac 1 4 \|[\X^\top \M \X, \N]\|_F^2.
\end{align*}
It follows from $\|\text{grad}~g(\X)\|_F^2 \leq \epsilon^2$ that
\begin{align}
\|(\X\X^\top - \I_N)\M\X \N\|_F^2 &\leq \epsilon^2, \label{eq:smallgrad1}\\
\|[\X^\top \M \X, \N]\|_F^2 &\leq 4 \epsilon^2. \label{eq:smallgrad2}	
\end{align}

Observe that
\begin{align*}
&\|[\X^\top \M \X, \N]\|_F^2 = \| \X^\top \M \X \N - \N \X^\top \M \X \|_F^2 \\
=& \sum_{i,j=1}^r (\X^\top \M \X)^2_{ij}	(\mu_j-\mu_i)^2
= \sum_{i,j=1 \atop i\neq j}^r (\X^\top \M \X)^2_{ij}	(\mu_j-\mu_i)^2 \\ \geq &\sum_{i,j=1 \atop i\neq j}^r (\X^\top \M \X)^2_{ij},
\end{align*}
where the above inequality follows from $(\mu_j-\mu_i)^2 \geq 1$ when $i\neq j$. Denote
\begin{align}
\X^\top \M \X = \D_{\M} + \E_{\M}	
\label{eq:xmxde}
\end{align}
with $\D_{\M}$ and $\E_{\M}$ being the diagonal and off-diagonal parts of $\M$, respectively. Then, we have
\begin{align}
\|\E_{\M}\|_F^2 \!=\!\! \sum_{i,j=1 \atop i\neq j}^r \!(\X^\top \M \X)^2_{ij} \!\leq \!\|[\X^\top \M \X, \N]\|_F^2 \!\leq\! 4\epsilon^2.	
\label{eq:embound}
\end{align}

Define $\E_1 \triangleq (\X\X^\top - \I_N)\M\X \N$. Then, we have
\begin{align}
(\X\X^\top - \I_N)\M\X = \X\X^\top \M\X - \M\X = \E_1\N^{-1},
\label{eq:xmxen}	
\end{align}
which immediately gives
\begin{align}
\!\!\!\|\X\X\!^\top\! \M\X \!-\! \M\X\|_F^2 \!=\! \!\|\E_1\N^{-1}\!\|_F^2 \!\leq \!\!\|\N^{-1}\!\|^2 \|\E_1\!\|_F^2 \!\leq \!\epsilon^2.	
\label{eq:xxmx}
\end{align}
On the other hand, by plugging $\X^\top \M \X = \D_{\M} + \E_{\M}$, we also have
\begin{align*}
&\|\X\X^\top \M\X - \M\X\|_F = \|\X\D_{\M} - \M\X + \X\E_{\M}\|_F \\
\geq &\|\X\D_{\M} - \M\X\|_F - \|\X\|\|\E_{\M}\|_F.	
\end{align*}
Combining this with \eqref{eq:xxmx} yields
\begin{align}
\|\X\D_{\M} \!-\! \M\X\|_F \!\leq \!	\|\X\X^\top \M\X\! -\! \M\X\|_F \!+\! \|\X\|\|\E_{\M}\|_F \!\leq \!3\epsilon,
\label{eq:xmxl1}
\end{align}
where we used $\|\X\| = 1$ and $\|\E_{\M}\|_F \leq 2\epsilon$.

Denote $\lambda_{i_j}$ as an eigenvalue of $\M$ that is closest to $\x_j^\top\M\x_j$, namely, $\lambda_{i_j} = \arg \min_{\lambda_n} (\x_j^\top \M \x_j - \lambda_n)^2$. Note that
\begin{equation}
\begin{aligned}
&\|\X\D_{\M} - \M\X\|_F ^2 = \sum_{j=1}^r \| \x_j (\x_j^\top \M\x_j) - \M\x_j \|_2^2	 \\
=& \sum_{j=1}^r \sum_{n=1}^N  ( \x_j^\top\M \x_j - \lambda_n )^2 \x_j(n)^2  \\
\geq & \sum_{j=1}^r \min_{1\leq n \leq N} ( \x_j^\top\M \x_j - \lambda_n )^2  \sum_{n=1}^N \x_j(n)^2 \\
=& \sum_{j=1}^r ( \x_j^\top\M \x_j - \lambda_{i_j} )^2,
\label{eq:xmxl2}
\end{aligned}
\end{equation}
where the last equality follows from the definition of $\lambda_{i_j}$ and $\|\x_j\|_2 = 1$. Denote $\bLambda_\Omega =  \diag([\lambda_{i_1},\cdots,\lambda_{i_r}])$ with $\lambda_{i_j} = \arg \min_{\lambda_n} (\x_j^\top \M \x_j - \lambda_n)^2$. Combining inequalities~\eqref{eq:xmxl1} and~\eqref{eq:xmxl2}, we can bound
\begin{equation}
\begin{aligned}
&\|\X^\top \M \X - \bLambda_\Omega\|_F^2 = \|\D_{\M} - \bLambda_\Omega\|_F^2 + \|\E_{\M}\|_F^2 \\
=& \sum_{j=1}^r ( \x_j^\top\M \x_j - \lambda_{i_j} )^2 + \|\E_{\M}\|_F^2
\leq  13\epsilon^2.
\label{eq:dmlo}
\end{aligned}
\end{equation}
Consequently, we can get~\eqref{eq:xMxeig1} by taking a square root on both sides of the above inequality.

It remains to show the inequality~\eqref{eq:dist_xv1}.
When $\Omega = [r]$, we have $\X_\Omega = [\I_r~ \zero]^\top \in\RRR^{N\times r}$. Denote $\X = [\X(\Omega,:)^\top~\X(\Omega^c,:)^\top]^\top$. Note that
\begin{align*}
\I_r = \X^\top \X = \X(\Omega,:)^\top \X(\Omega,:) + \X(\Omega^c,:)^\top \X(\Omega^c,:).
\end{align*}
Taking the trace of both sides yields
\begin{align*}
r &= \tr(\I_r) = \tr\left(\X(\Omega,:)^\top \X(\Omega,:)\right)	+ \tr\left(\X(\Omega^c,:)^\top \X(\Omega^c,:)\right)\\
&=\| \X(\Omega,:) \|_F^2 + \| \X(\Omega^c,:) \|_F^2\\
&= \sum_{j=1}^r \X^2(j,j) + \|\X_\Omega^{\text{off}}\|_F^2 + \| \X(\Omega^c,:) \|_F^2,
\end{align*}
where $\X_\Omega^{\text{off}}\in\RRR^{r\times r}$ denotes the off-diagonal part of $\X(\Omega,:)\in \RRR^{r \times r}$. It follows that
\begin{equation}
\begin{aligned}
&\|\X_\Omega^{\text{off}}\|_F^2 + \| \X(\Omega^c,:) \|_F^2
=  r - \sum_{j=1}^r \X^2(j,j) \\
 	= & \sum_{j=1}^r\!\! \left( 1 \!+ \!\X^2(j,j) \!-\! 2\X^2(j,j)  \right)
\geq \! \sum_{j=1}^r \!\left( 1  \!-\! \X(j,j)  \right)^2,
\label{eq:xooff}
\end{aligned}
\end{equation}
where the last inequality follows from $0\leq \X(j,j) \leq 1$.\footnote{Note that we assume $\X(j,j) \geq 0$. In the case when $\X(j,j) < 0$, one can let $\X_\Omega(j,j) = -1$.}

Note that
\begin{align*}
\|\X\D_{\M} - \M\X\|_F =& \|\X\bLambda_{\Omega} - \M\X + \X(\D_{\M} - \bLambda_{\Omega} )\|_F \\
\geq & \|\X\bLambda_{\Omega} - \M\X\|_F - \|\X(\D_{\M} - \bLambda_{\Omega} )\|_F,
\end{align*}
which together with~\eqref{eq:xmxl1} gives us
\begin{equation}
\begin{aligned}
\|\X\bLambda_{\Omega} \!-\! \M\X\|_F \!\leq &  \|\X\D_{\M} \!-\! \M\X\|_F \!+\! \|\X(\D_{\M} \!-\! \bLambda_{\Omega} )\|_F \\
\leq & 3\epsilon + 	\|\D_{\M} - \bLambda_{\Omega} \|_F \leq 6\epsilon,
\label{eq:xlomx}
\end{aligned}
\end{equation}
where the last inequality follows from combining~\eqref{eq:xmxl2} and~\eqref{eq:xmxl1}.
It follows from~\eqref{eq:xlomx} that
\begin{align*}
\|\X\bLambda_{\Omega} - \M\X\|_F \!=\! & \left\|\! \left[\!
\begin{array}{c}
\X(\Omega,:) \bLambda_\Omega - \bLambda_\Omega \X(\Omega,:)\\
\X(\Omega^c,:) \bLambda_\Omega - \bLambda_{\Omega^c} \X(\Omega^c,:)
\end{array}
\!\right]\!\right\|_F \leq 6\epsilon,
\end{align*}
where $\bLambda_{\Omega^c}\in \RRR^{(N-r)\times (N-r)}$ is a diagonal matrix that contains the remaining $N-r$ eigenvalues of $\M$. Then, we have
\begin{align*}
36\epsilon^2 \geq &\| \X(\Omega,:) \bLambda_\Omega \!-\! \bLambda_\Omega \X(\Omega,:)  \|_F^2	\!=\! \| \X(\Omega,:) \odot [\lambda_j - \lambda_i] \|_F^2 \\
\geq &d_{\min}^2 \|\X_\Omega^{\text{off}}\|_F^2,\\
36\epsilon^2 \geq & \| \X(\Omega^c\!,:) \bLambda_\Omega \!- \!\bLambda_{\Omega^c} \X(\Omega^c\!,:)  \|_F^2	\!=\! \| \X(\Omega^c\!,:) \!\odot \![\lambda_j \!-\! \lambda_i^c] \|_F^2 \\\geq & d_{\min}^2 \|\X(\Omega^c,:)\|_F^2,
\end{align*}
where $[\lambda_j - \lambda_i]\in\RRR^{r\times r}$ and $[\lambda_j - \lambda_i^c]\in \RRR^{(N-r)\times r}$ denote two matrices with the $(i,j)$-th entry being $\lambda_j - \lambda_i$ and $\lambda_j - \lambda_i^c$, respectively. Here, we use $\lambda_i$ and $\lambda_i^c$ to denote the $i$-th diagonal entry of $\bLambda_\Omega$ and $\bLambda_{\Omega^c}$, respectively. Adding the above two inequalities gives
$\|\X_\Omega^{\text{off}}\|_F^2 + \|\X(\Omega^c,:)\|_F^2 \leq 72 d_{\min}^{-2} \epsilon^2,$
which together with~\eqref{eq:xooff} gives
\begin{align*}
\|\X-\X_\Omega\|_F^2 =& 	\|\X(\Omega^c,:)\|_F^2 + \|\X_\Omega^{\text{off}}\|_F^2 + \sum_{j=1}^r \left( 1  - \X(j,j)  \right)^2\\
\leq &2(\|\X(\Omega^c,:)\|_F^2 + \|\X_\Omega^{\text{off}}\|_F^2) \leq 144 d_{\min}^{-2} \epsilon^2.
\end{align*}

\end{document}